\documentclass[12pt]{article}
\usepackage[utf8]{inputenc}
\usepackage[T1]{fontenc}
\usepackage{amsfonts,amssymb,amsmath}
\usepackage{geometry}
\usepackage{graphicx}
\usepackage{titling}
\usepackage{amsthm}
\usepackage{bbm}
\usepackage{mathtools}
\usepackage{bm}

\usepackage{pdfpages}

\usepackage[colorlinks=true, citecolor=blue]{hyperref}
\hypersetup{
    colorlinks=true,
    linkcolor=blue,
    filecolor=blue,      
    urlcolor=blue,
}

\usepackage{dsfont}

\usepackage{enumitem}
\setlist[itemize,1]{label=\textendash}

\def\P{{\mathbb P}}
\def\N{{\mathbb N}}
\def\Z{{\mathbb Z}}
\def\R{{\mathbb R}}
\def\E{{\mathbb E}}

\newcommand{\p}{\textup{\texttt{+}}}
\newcommand{\m}{\textup{\texttt{-}}}

\newcommand{\term}[1]{\textbf{#1}}   
\newcommand{\foreign}[1]{\emph{#1}}  

\usepackage[giveninits=true]{biblatex}
\DeclareNameAlias{sortname}{last-first}

\title{The Pathwise Approach to Metastability and its Applications to Galves--Löcherbach Models}
\date{\today}
\author{Morgan André\thanks{Departamento de Estatística, Instituto de Matem\'atica e Estat\'istica, Universidade de São Paulo, Rua do Matão 1010, 05508-090, Brazil. E-mail: morgan.andre@usp.br} \ and Kádmo Laxa\thanks{Faculdade de Filosofia, Ciências e Letras de Ribeirão Preto, Universidade de São Paulo, Av. Bandeirantes, 3900, Ribeirão Preto-SP, 14040-901, Brazil. E-mail: kadmo.laxa@usp.br}}

\begin{document}

\maketitle

\begin{abstract}
Metastability is the tendency of a system to dwell for a very long time near an apparently stable equilibrium before a rare fluctuation drives it, on a comparatively short time scale, towards another. Among the rigorous frameworks developed to capture this phenomenon, the pathwise approach proceeds by identifying the ``typical'' trajectories of the stochastic dynamics at hand and estimating their probabilities. In this article we review the pathwise approach and its application to the Galves--Löcherbach (GL) class of stochastic models of spiking neural networks. After recalling the conceptual and historical roots of the theory --- which goes from chemistry to rigorous probability theory, with fundamental ideas coming mainly from statistical physics --- and illustrating them on two classical examples, we give a general definition encompassing the known variants of the GL model and survey the metastability results already established for some of these variants. As far as we can, we do so in a self-contained fashion, and we sketch the proofs when possible, highlighting their common structure. We close with a discussion on open problems and point to possible further directions.
\end{abstract}

\section{Introduction}

The phenomenon of metastability has been widely observed throughout history across many distinct fields, originating primarily in physics and chemistry. Phenomenologically, it may be roughly described as the tendency of a system --- a solution, a fluid, a ferromagnetic material, and so on --- to remain for a long time in an apparently stable equilibrium, until the occurrence of an unlikely perturbation triggers a transition toward another equilibrium --- typically stable, though not necessarily --- on a much shorter time scale. In most cases the perturbation is endogenous to the dynamics of the system, giving the observer the puzzling impression that the sudden transition occurred spontaneously.

One of the simplest and most well-known examples of such behavior is that of supercooled water, which arises when the temperature of a bottle of distilled water is brought slightly below $0\,^{\circ}$C. The purity of the water permits the persistence of a metastable liquid state, even though one would expect a transition toward the solid state, which is the only stable equilibrium. Nonetheless, the liquid state stubbornly persists for a possibly long time, until the transition occurs \emph{suddenly} at a seemingly random moment.

Another historically important manifestation of metastability arises in chemistry when two reactive compounds are mixed. The reaction may not occur immediately, but instead be delayed for an apparently unpredictable duration. This delay can be explained as follows: for the reaction to proceed, colliding molecules must possess sufficient kinetic energy, and the proportion of such molecules in the solution might be small, so that one has to wait for this random event to happen. In other words, the reaction rate can be characterized by the statistical properties of the solution, viewed as a population of molecules randomly colliding with one another. This delay has been characterized quantitatively by the Arrhenius law, which nonetheless initially emerged as an essentially empirical relation, and a mathematically rigorous explanation of such phenomena long remained out of reach for the deterministic framework of the natural sciences. Such an explanation could not arise before the revolutionary work of Ludwig Boltzmann on the kinetic theory of gases had percolated deeply enough into the modern scientific mindset, opening physics to the --- still not fully rigorous --- ideas of statistics and probability theory, and laying the groundwork for statistical mechanics. Einstein's explanation of Brownian motion in 1905 was a particularly important milestone along this path, which led, a few decades later, to the first successful attempt at a rigorous stochastic justification of the metastable reactions observed in chemistry. This was achieved through a simple \emph{mesoscopic} model --- a Brownian particle in a double-well potential --- introduced by Hendrik Kramers in 1940. This quantitative framework provided a fundamental insight into the nature of metastable macrostates: they can essentially be understood as local minima of the energy landscape along which the system evolves. Whenever one of these local minima is reached, the system naturally tends to remain there for a long time, until a fluctuation --- unlikely at a microscopic time scale but unavoidable at a macroscopic one --- causes the crossing of the energy barrier separating this local minimum from another. Kramers' work eventually led to the development of the theory of random perturbations of dynamical systems by Freidlin and Wentzell, inaugurating the pathwise approach to metastability, in which the most likely trajectories of the stochastic dynamical system are identified and their probabilities estimated.

In parallel, efforts were made in the field of interacting particle systems to propose a rigorous framework for understanding metastability in models describing complex systems at the microscopic scale, ultimately leading to the characterization proposed by Penrose and Lebowitz \cite{evolutionofensembles} in the early 1970s. This characterization decomposed the phenomenon of metastability into three distinct criteria:
\begin{enumerate}
\item the existence of a unique stable state,
\item the tendency of the system to remain for a very long time outside this equilibrium,
\item and the ``irreversibility'' of the transition from the metastable to the stable state.
\end{enumerate}

The approach developed by Penrose and Lebowitz --- known as the evolution of ensembles --- did not however fully meet the axiomatic standards of mathematical theories, even though the field of Markovian interacting particle systems (IPS) had already been established as a branch of probability theory in the late 1960s, following independent foundational contributions by F.\ Spitzer and R.\ Dobrushin. A mathematically rigorous treatment was ultimately achieved in 1984, when the pathwise approach of Freidlin and Wentzell and the ideas of Penrose and Lebowitz were combined into a rigorous characterization of metastability for IPS in the seminal work of Cassandro, Galves, Olivieri, and Vares \cite{cassandro}. This work illustrated the approach through the rigorous study of the metastable properties of two important models: the Curie--Weiss model and the contact process, and served as a point of departure for the subsequent analysis of metastability in various other classical complex systems.

Metastability has since been observed and studied across a wide range of fields, including economics, electronics, and population biology, among others --- spanning the entire spectrum of rigor, from qualitative descriptions to quantitative mathematical frameworks. In recent years, metastability has also emerged as a prominent concept in neuroscience, and it is now widely accepted that it is one of the brain's fundamental dynamical properties. The precise meaning of this claim, however, remains to be fully clarified. A key observation is that --- in order to adapt consistently to noisy and ever-changing environments --- brain sub-networks exhibit two competing tendencies: on the one hand, they express independently specialized functions at mesoscopic time scales (segregation), and, on the other hand, they cooperate at macroscopic time scales (integration). Borrowing terminology from statistical physics, one could reformulate this by saying that brain dynamics must strike a balance between stability and instability. This delicate equilibrium is, broadly speaking, what is referred to as metastability in the brain.

The problem of rigorously characterizing metastability in that context has been approached in a variety of ways, most notably through dynamical systems theory (see for example \cite{rabinovich2008transient} or \cite{beim2019metastable}) and coordination dynamics (see for example \cite{metabrain}). These approaches have, however, remained somewhat loosely connected to the well-established notion of metastability formalized in statistical physics, which carries the considerable advantage of being underpinned by a powerful mathematical framework. Meanwhile, evidence has accumulated over the years that neural systems exhibit intrinsically stochastic dynamics --- both at the biophysical level, through observations of the random activation of ion channels \cite{colquhoun}, and at the cognitive level, through arguments on the nature of learning \cite{deco}. To account for this stochasticity, \textcite{glmodel} introduced a probabilistic model of spiking neurons (see also the review \cite{glarticle}), which has subsequently been studied in multiple variants\footnote{The original model is discrete-time, but subsequent variants have been studied in a continuous-time framework; some are formulated as particle systems, others as piecewise-deterministic Markov processes, and so on. See \cite{gl3,gl4,evafour,evamonm,baccelli1,baccelli2,mandre1,mandre2,mandre3,ferrari,karina,lud,romaro,laxa,pouzatandre,duarteost,amarcos}.}, and has become generically known as the Galves--Löcherbach (GL) model. The continuous-time variant of this model can essentially be seen as a system of interacting point processes similar to multivariate Hawkes processes, with an additional post-spiking reset mechanism.  Notably, one of the defining features of this class of models is that it is amenable to the rigorous approaches to metastability developed in the statistical physics literature, and it has accordingly been studied from this perspective. The main objective of this work is to provide a concise yet reasonably comprehensive account of these developments, thereby establishing the state of the art of this line of research.

This article is organized as follows. In Section \ref{sec:pathwise} we give a short historical review of the mathematical approaches leading to the pathwise approach to metastability and state the rigorous definition of metastability proposed by this approach. In Section \ref{sec:GLModels}, after a short biologically motivated introduction to neuronal networks, we give a general definition of the Galves-Löcherbach framework. In Section \ref{sec:metastability_GLmodel} we present the main metastability results for GL models available in the literature and discuss the common ideas and strategies mobilized for proving these results. In Section \ref{sec:conclusion} we present some open problems and make some final observations.

\section{The pathwise approach to metastability} \label{sec:pathwise}

In this section, we give a short historical review of the mathematical approaches that led to the pathwise approach to metastability, before stating a reasonably general and rigorous definition of metastability in Section~\ref{subsec:pathwisedef}.

\subsection{Evolution of the ensembles}

In \cite{evolutionofensembles}, metastability is studied from the ``evolution of ensembles'' point of view. The authors introduce a general method for describing metastable states in statistical mechanics and use this method to study metastability for the generalized van der Waals system describing a system of interacting particles.

Consider a thermodynamic system given by an initial configuration and a law describing the evolution. The initial configuration is given by a probability measure in the space of all possible configurations, and the evolution of the system can be deterministic (dynamical system) or stochastic (Markov process). The generalized van der Waals system considered in \cite{evolutionofensembles} is an example of a deterministic evolution, while  \textcite{evolutionofensembles2} considers a Markovian evolution for a two-dimensional Ising spin system.

In order to define metastability in this context one needs to define a suitable Gibbs measure $\mu_R$ with support in a suitably chosen region $R$ in the space of all configurations. The region $R$ and the measure $\mu_R$ are said to be metastable for a given thermodynamic system if the following three properties are satisfied.
\begin{enumerate}
\item Only one thermodynamic phase is present.
\item A system with initial distribution $\mu_R$ is likely to take a long time to get out of $R$.
\item Once the system has gotten out, it is unlikely to return to $R$.
\end{enumerate}

In order to deal with the second point above, \cite{evolutionofensembles} consider the probability $p_t$ of being outside $R$ at time $t$ with initial configuration $\mu_R$. For the models considered in \cite{evolutionofensembles} and \cite{evolutionofensembles2} it is shown that the maximal increase rate of $p_t$ is achieved at $t=0$. This value is called the \term{escape rate} and is pivotal to their approach. Formally, denoting this escape rate by $\lambda$, one has:
$$
\lambda \triangleq \sup_{t \in \R_+} \frac{dp_t}{dt} = \left. \frac{dp_t}{dt} \right|_{t=0}.
$$ 

The second property in the above characterization of metastability can be made precise by requiring that the escape rate $\lambda$ becomes arbitrarily small for an appropriate choice of the model parameters (such as the number of particles, the size of the configuration space, and $R$).

To deal with the third property, \cite{evolutionofensembles} argue that the  measure of the configurations of $R$ must be negligible at equilibrium. In both the model with deterministic evolution studied in \cite{evolutionofensembles} and in the Markovian system considered in \cite{evolutionofensembles2}, the measure of the configurations on $R$ under the Gibbs measure taken over the whole set of configurations is negligible. Formally, for a suitable choice of parameters of the model (the same parameter considered above which makes the escape rate as small as we want), we can make this measure as small as we want.  Note that in the Markovian system the Gibbs measure is the unique invariant measure of the considered Markov process. 

From the ``evolution of ensembles'' point of view, the study of metastability concerns the study of the equilibrium properties of the system and the time evolution of the probability distribution in the space of all possible configurations, given an initial probability distribution restricted to a suitable region. Therefore, the evolution of a typical trajectory of the system is not analyzed.

In \cite{cassandro} Cassandro et al.\@ observe that, by analyzing only the time evolution of the \emph{probability distribution}, one cannot distinguish between the genuine metastable case in which single trajectories tend to stay in a ``stable'' situation for a long and unpredictable time before making a sharp transition to the true equilibrium, and the case in which single trajectories present a smooth, very slow evolution toward the stable situation, which is completely unrelated to metastability.  The pathwise approach to metastability introduced by \cite{cassandro} overcomes this problem by identifying and studying the typical trajectories of the process.

\subsection{Freidlin--Wentzell theory}

The Freidlin--Wentzell theory (see \cite{fwbook} for a complete account of this theory) deals  with Markov processes defined by the following equation
\begin{equation}\label{eq:FWprocess}
X^{\epsilon}(t)=x+\int_0^tb(X^{\epsilon}(s))ds+\epsilon W(t), \quad t\geq0,  
\end{equation}
where $\epsilon>0$, the initial position is $x\in \R^d$, $d\geq 1$, $b:\R^d\to \R^d$ is a globally Lipschitz function and $(W(t))_{t \in \R_+}$ is a standard d-dimensional Brownian motion starting at the origin.  Informally, $(X^{\epsilon}(t))_{t \in \R_+}$ is a small random perturbation of the solution of the deterministic equation 
\begin{equation}\label{eq:deterministic}
\frac{dx(t)}{dt} =b(x(t)) \quad \text{ and } \quad x(0) =x.
\end{equation}
The random perturbation is given by $\epsilon W(t)$, which means that it vanishes as $\epsilon \to 0$. With this, we expect that for any fixed $T>0$, $(X^{\epsilon}(t))_{t\in [0,T]}$ is close to $(x(t))_{t\in [0,T]}$ if $\epsilon$ is small. In fact, for any $\delta>0$, there exist $C_1,C_2>0$ such that
$$
\P\left(\sup_{t\in [0,T]}|X^{\epsilon}(t)-x(t)|>\delta\right) \leq C_1\delta^{-1}e^{-C_2\delta^2\epsilon^{-2}} \to 0, \text{ as } \epsilon\to 0.
$$
The constants $C_1$ and $C_2$ depend only on $d,T$ and the Lipschitz constant of $b$ (see equation (2.35) of \cite{metaest2}). Since in the inequality above the right-hand side does not depend on the initial position $x$, we can take the supremum for $x\in \R^d$ on the left-hand side. 

\begin{figure}[ht!]
\centering 
\includegraphics[width=0.7\textwidth]{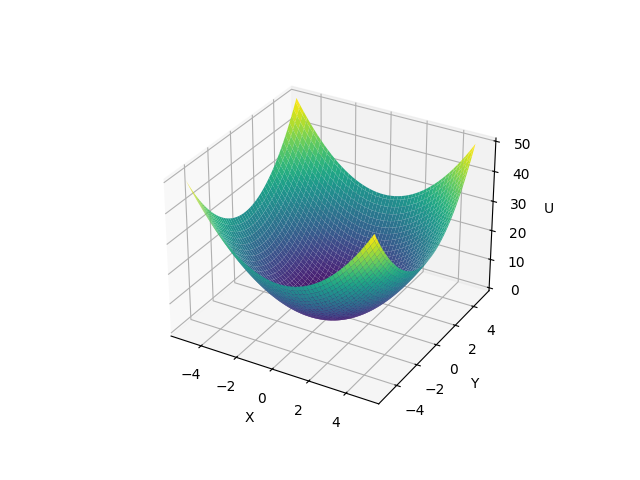}\\
\caption{\label{fig:x2y2} Plot of a function $U:\R^2\to \R$. In the deterministic equation \eqref{eq:deterministic} with $b(x)=-\nabla U(x)$, all points around the origin are attracted to the origin.}
\end{figure}

In the following, let us consider the case $b(x)=-\nabla U(x)$, where $U:\R^d\to \R$ is a function of class $C^2$. In this case, the plot of the function $U$ illustrates the behavior of the deterministic equation \eqref{eq:deterministic}: it follows the direction in which $U$ decreases the most. Figure \ref{fig:x2y2} illustrates the case in which all points around the origin are attracted to the origin.

The theory developed by Freidlin and Wentzell gives a detailed description of how and when $(X^{\epsilon}(t))_{t \in \R_+}$ exits an open and bounded domain $D\subset \R^d$ attracted to a point $p \in D$. Consider that for any initial point in $D$, the deterministic process $x(t)$ approaches $p$ as $t\to \infty$. Moreover, $(x(t))_{t \in \R_+}$ never exits $D$.
Using large deviation estimates (see Theorem 2.25 and Definition 2.31 in \cite{metaest2}), one can introduce for $x,y \in D$ a non-negative value $V(x,y)$ that is the ``cost'' that represents how much the process \eqref{eq:FWprocess} with initial point $x$ has to move away from its limit behavior \eqref{eq:deterministic} to be close to $y$. Under additional assumptions on the trajectories $(x(t))_{t \in \R_+}$, the set $D$ and its boundary $\partial D$ (see subsection 2.6 of \cite{metaest2}), Freidlin and Wentzell prove the following.

Consider an initial point $x\in D$ and assume that $y \in \partial D$ is the unique point that minimizes $V(x,z)$ for all $z\in \partial D$. The point of $\partial D$ at which $(X^{\epsilon}(t))_{t \in \R_+}$ exits $D$ for the first time gets closer to $y$ as $\epsilon\to 0$. 
Formally, for any $\delta>0$,
denoting by $\tau_{\epsilon}$ the first time at which $(X^{\epsilon}(t))_{t \in \R_+}$ exits $D$, we have that
$$
\lim_{\epsilon\to 0}\P(|X^{\epsilon}(\tau_{\epsilon})-y|>\delta)=0.
$$
Moreover, by denoting $V_0\triangleq\inf\{V(x,z):z\in \partial D\}$, we have that
$$
\lim_{\epsilon\to 0} \epsilon^{2} \log \E(\tau_{\epsilon})=V_0,
$$
and
$$
\lim_{\epsilon\to 0}\P(e^{\epsilon^{-2}(V_0-\delta)}< \tau_{\epsilon} <e^{\epsilon^{-2}(V_0+\delta)})=1, \text{ for any } \delta>0.
$$
With this we show that $\E(\tau_{\epsilon})$ is exponentially large with order $V_0$ as $\epsilon\to 0$.

\begin{figure}[ht!]
\centering 
\includegraphics[width=0.7\textwidth]{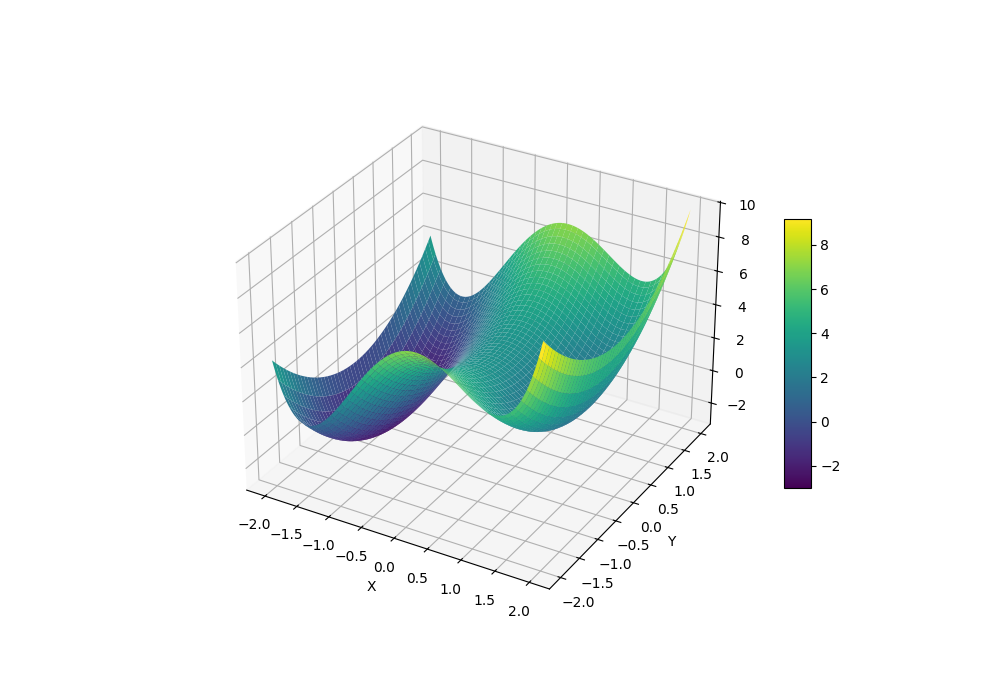}\\
\caption{\label{fig:doublewell3d} Plot of a function $U:\R^2\to \R$. In the deterministic process \eqref{eq:deterministic} with $b(x)=-\nabla U(x)$, all points around the bottom of each well are attracted to the bottom of the respective well.}
\end{figure}

The study of metastability is considered when we have a double-well structure: two disjoint sets, each one attracted to a different stable point. This situation is illustrated by Figure \ref{fig:doublewell3d}. If the process starts close to the bottom of the higher well, the system exhibits a metastable behavior for $\epsilon>0$ small, remaining a long time in the higher well before making a transition to the set of points attracted to the deeper well. The results above, giving us a precise description of how and when the process exits an open domain, are the tools to describe this metastable behavior by studying the most likely trajectories of the process  (see Chapter 5 of \cite{metaest2}).

\subsection{The pathwise approach to metastability} \label{subsec:pathwisedef}

In this section we give a reasonably general and rigorous definition of metastability, inspired by the framework introduced in \cite{cassandro} as part of the pathwise approach. This perspective is the one to be remembered for the following sections, since it is the framework in which most metastable results on GL models were obtained. Consider a stochastic process $(X^N(t))_{t \in \R_+}$ taking values in some space $\mathcal{S}_N$, where $N$ represents the \term{volume} of the system --- i.e.\@ the number of particles, neurons, individuals etc. Quite often it will be useful to think of $\mathcal{S}_N$ as \emph{embedded} in a wider metrizable space $\mathcal{S}$. The precise meaning of the embedding depends on the context, but in most cases $\mathcal{S}$ is either an infinite cartesian product (or equivalently a power set) --- and in that case $\mathcal{S}_N$ is the projection of $\mathcal{S}$ on finitely many coordinates --- or simply a (possibly uncountable) set of numbers --- and in that case $\mathcal{S}_N$ is a finite subset of $\mathcal{S}$. As it will become clear below, metastability is an asymptotic property in the $N \rightarrow \infty$ limit, which is why the limit space $\mathcal{S}$ often plays an important role. In other contexts metastability can be expressed relative to another parameter, different from the volume, such as inverse temperature, viscosity, noise intensity and so on. Through this review the relevant parameter will almost always be the volume $N$, with some exceptions in Section \ref{sec:related}.

One needs to identify two disjoint subsets $\mathcal{M}_N \subset \mathcal{S}_N$ and $\mathcal{A}_N \subset \mathcal{S}_N$, representing respectively the \term{metastable region} and the \term{stable region}, or ``absorbing region''. $\mathcal{M}_N$ and $\mathcal{A}_N$ are frequently complementary sets to each other, but not always. In some cases it might be useful to consider the corresponding limit sets $\mathcal{M} \subset \mathcal{S}$ and $\mathcal{A} \subset \mathcal{S}$, whatever the ``limit'' formally means in that context --- typically either the projective limit, the Hausdorff limit or simply the union over all $N \in \N$. The process starts from some fixed initial state $x \in \mathcal{M}_N$, and stays there during thermalization, behaving in a way \emph{which resembles stationarity}. At some point nonetheless, a tunneling occurs, that is: \emph{a rapid and unexpected transition} from the metastable region $\mathcal{M}_N$ toward the stable region $\mathcal{A}_N$. Let $\tau_N$ denote the tunneling time, that is, the hitting time of $\mathcal{A}_N$ by $(X^N(t))_{t \in \R_+}$. Below any statement is implicitly expected to hold with respect to the probability measure $\P=\P_x$, under which $X^N(0)=x$ almost surely, where $x$ is any given state in $\mathcal{M}_N$. The  characterization of metastability, as proposed by the pathwise approach, is aimed at making the two claims above precise and rigorous. Moreover one has to identify the relevant \emph{time-scales} on which these phenomena occur. Besides the standard time-scale $t \geq 0$ by which the system is indexed --- called the \term{microscopic} time-scale --- there must be a macroscopic one, on which one might observe the tunneling outside the metastable region. Furthermore there must be a third, intermediary time-scale --- the mesoscopic one --- on which one observes thermalization. Formally, metastability is characterized by the two following properties.

\begin{enumerate}[label=\textbf{P.\arabic*}]
    \item \label{meta:loss} \textbf{Asymptotic loss of memory:} There exists a \term{macroscopic} time scale $(\beta_N)_{N \in \N}$ on which $(\tau_N/\beta_N)_{N \in \N}$ converges in distribution to a unit mean exponential random variable. That is, for all $t > 0$, one has:
    $$ \lim_{N \to \infty}\P\left(\frac{\tau_N}{\beta_N}>t\right)=e^{-t}. $$

    \item \label{meta:therma} \textbf{Thermalization:} There exists a probability measure $\mu$ on $\mathcal{S}$ and a \term{mesoscopic} time scale $(R_N)_{N \in \N}$, diverging as $N \to \infty$, but not as fast as the macroscopic time scale --- that is $R_N = o(\beta_N)$ --- such that the system is approximately described on this time-scale by the restriction of $\mu$ to $\mathcal{S}_N$ --- which puts mass $0$ on $\mathcal{A}_N$. By this we mean that for any $f$ in a wisely chosen class of real functions on $\mathcal{S}$, the \emph{temporal averages} with respect to $f$ of $(X^N(s))_{s \in \R_+}$ on this mesoscopic time-scale converge in probability to the \emph{spatial average} of $\mu$ with respect to $f$. More precisely, for suitable $f$ and arbitrarily small $\epsilon > 0$, define:
    \begin{equation} \label{eq:typical}
    \Gamma_N^ \epsilon(f) \triangleq \left\{\max_{\substack{t<\tau_N-R_N \\ t \equiv 0 \, \text{mod} \, R_N}}\left|\frac{1}{R_N}\int_{t}^{t+R_N} f\left(X^N(s)\right) \, ds - \int f \, d\mu\right| < \epsilon \right\}.
    \end{equation}
    When $\mathcal{S}$ is an infinite cartesian product we require that the function $f$ be \emph{cylindrical}\footnote{A cylindrical function is one whose value depends only on finitely many coordinates.}, to ensure that $\Gamma_N^ \epsilon(f)$ is well defined at least for big enough $N$. Then one has:
    $$ \lim_{N \to \infty} \P\left(\Gamma_N^ \epsilon(f) \right) = 1. $$
\end{enumerate}

Property \ref{meta:loss} captures the abruptness in the transition out of the metastable phase, since memorylessness is the probabilistic translation of unpredictability: for $\tau_N$ to have exponential distribution means that knowing that the system is still in its pseudo-stationary phase after a given amount of time tells you nothing about how much more you shall wait before seeing the transition. While it might not be obvious at first sight, the ``pathwiseness'' already plays a role in the memorylessness property: hitting times in general can indeed be seen as functionals over the trajectories of the system, so that the distribution of $\tau_N$ shall yet encode pathwise information not available at the level of the time-dependent probability distribution of the evolution of ensembles approach. Property \ref{meta:therma} on the other hand relates to single trajectories quite straightforwardly. The set  $\Gamma_N^ \epsilon (f)$ --- namely, the set of trajectories whose temporal averages with respect to $f$ stay close to the prescribed value $\int f \, d\mu$ along time --- defines the typical trajectories. The function $f$ is taken out of a specific class of functions that need to be specified --- and it represents what one would call an observable quantity in physics context. Moreover, the identification of the typical trajectories essentially boils down to the identification of the pseudo-stationary measure  $\mu$ from which $\Gamma_N^ \epsilon (f)$ is defined. Finally, one has to justify the adjective ``typical'' by showing that the trajectories outside this class appear with vanishingly small probability.

While the exact formulation of \ref{meta:therma} presented above is somewhat technical, it shall not hide the simple fundamental idea, which is the existence of a probability measure $\mu$ giving an approximate description of the system along time prior to absorption. The fact that the description is approximate is what captures the pseudo-stationarity of the dynamics --- since if the description were exact then $\mu$ would simply be an invariant measure for the process under consideration. In many of the later results presented in this review (those reviewed in section \ref{sec:metastability_GLmodel} in particular), while \ref{meta:therma} is usually not formally proven in the sense prescribed above, the fundamental idea of having an invariant measure $\mu$ giving an approximate description of the system prior to extinction is almost always present.

Properties \ref{meta:loss} and \ref{meta:therma} can be neatly mixed together into a single statement capturing the essence of metastability. For any positive real numbers $R$ and $t$, we define the spatio-temporal average of $(X^N(t))_{t \in \R_+}$ on the time interval $[t, t+R]$ as a random probability measure on $\mathcal{S}_N$ given by:
\begin{equation} \label{eq:stavg}
A_R^N (t, \cdot) \triangleq \frac{1}{R} \int_t^{t+R} \delta_{X^N(s)} (\, \cdot \,) ds.
\end{equation}
Above $\delta_x(\, \cdot \,) = \mathbbm{1}_{x \in \, \cdot \,}$ denotes the Dirac measure on $x$. Then $(A_R^N(t))_{t \in \R_+}$ can be seen as a stochastic process itself, taking value in the space of probability measures on $\mathcal{S}$. Then, giving some topological structure\footnote{Like the $w^*$-topology.} to this space, one can usually show by standard tightness arguments that \ref{meta:loss} and \ref{meta:therma} imply convergence of $(A_{R_N}^N(\beta_Nt))_{t \in \R_+}$ toward a given limit $(A(s))_{s \in \R_+}$ in the usual Skorohod space, typically of the following form: $$A(s) = \begin{cases}
\mu &\text{ if } s < T,\\
\nu &\text{ if } s \ge T.\\
\end{cases}$$
Above $T$ is an exponential random variable of rate $1$, and $\nu \neq \mu$ is some probability measure on $\mathcal{A}$ describing the true equilibrium of the system. The two measures $\mu$ and $\nu$ shall hence be interpreted as two \emph{macrostates}, visible only at the macroscopic scale on both space and time, and describing two qualitatively different phases of the system, visited in sequence. While $\mu$ represents the metastable phase, i.e.\@ some false equilibrium, $\nu$ represents the only legitimate equilibrium of the system.

\subsection{Two introductory examples}

The very first models that have been shown to present metastable behavior in the pathwise sense are the Curie--Weiss model and the classical Contact process. The main part of the analysis was done in \cite{cassandro}. In this section we introduce these two models and briefly highlight the results obtained and techniques involved in this seminal paper. For the sake of conciseness the ideas behind the proofs are mentioned quite allusively below, and the proofs themselves are not really discussed at all. We refer to chapter 4 of \cite{metaest2} for a modern in-depth exposition.

\subsubsection*{The Curie--Weiss model}
The Curie--Weiss model is a \emph{mean-field} simplification of the Ising model. Consider a system of $N$ particles with spins $\sigma = (\sigma_1, \ldots, \sigma_N) \in \{\m 1,\p 1\}^N$, whose interactions are described by the following Hamiltonian $$H(\sigma) \triangleq -\frac{1}{N} \sum_{i<j} \sigma_i \sigma_j - h \sum_i \sigma_i - \frac{1}{2}.$$
The difference between the Curie--Weiss model and the classical Ising model lies in the first term, which represents the interactions, and in that case is an average over all spins in the system, while the average would be taken only over immediate neighbors in the classical Ising model. This means we are assuming complete interaction. Then the Gibbs measure for temperature $T$ is given by $ \mu(\sigma) \propto e^{-H(\sigma)/T}$. Since the system has complete interaction, for a given number $n$ of particles with positive spins, the exact configuration $\sigma$ is unimportant: $H(\sigma) = - \tfrac{1}{2N}(2n - N)^2 - (2n - N)h \triangleq U(n)$.
Therefore, the Gibbs measure $\mu$ can equivalently be seen as a measure on the space $\{0,1, \ldots, N\}$ of the number of positive spins: $\mu (n) \propto \binom{N}{n} e^{- U(n)/T}$.

\begin{figure}[ht!]
\centering 
\includegraphics[width=0.5\textwidth]{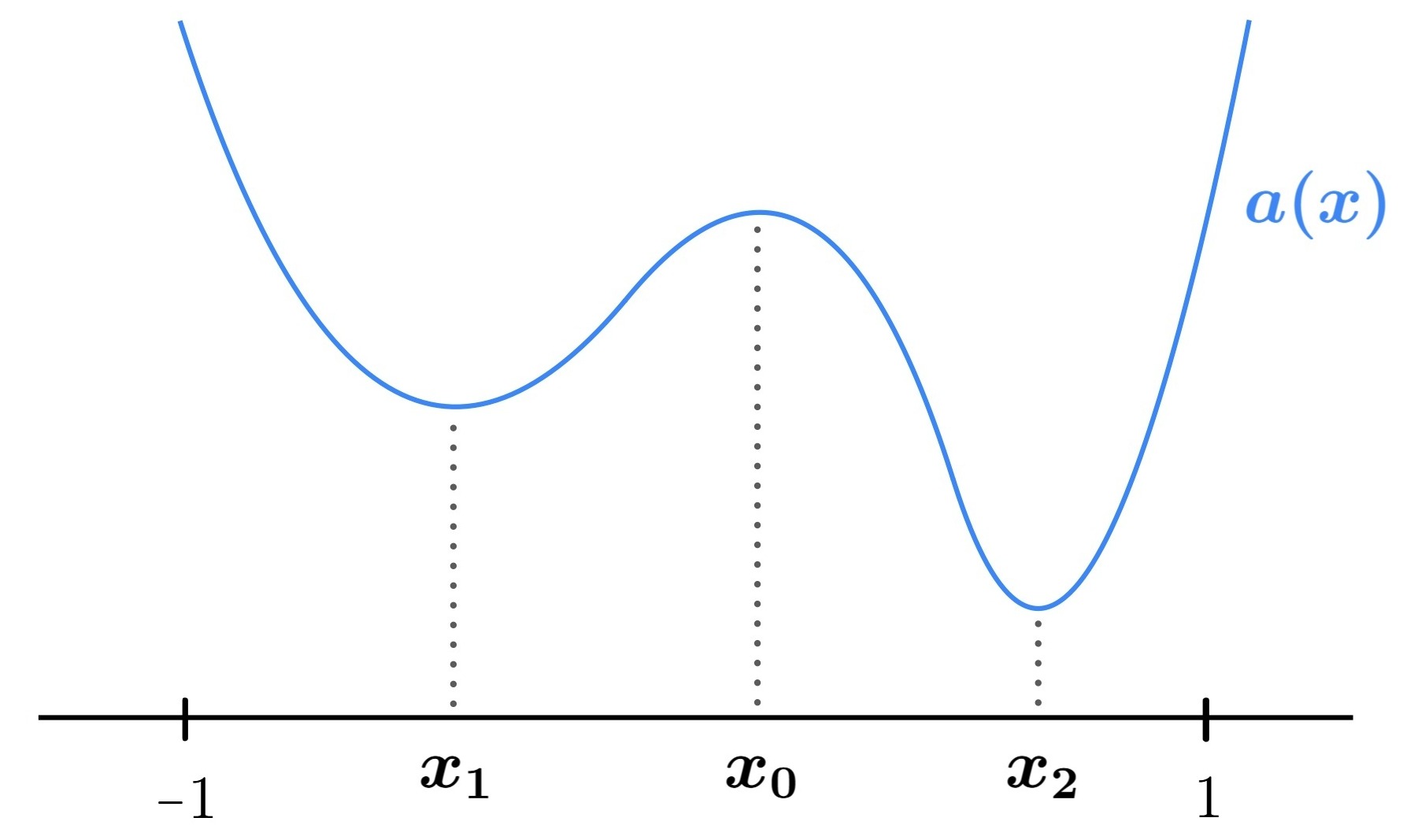}\\
\caption{\label{fig:double-well} The energy profile $a(x)$ as a function of the rescaled magnetization $x$. $x_1$ and $x_2$ denotes the position of the two wells, while $x_0$ is the position of the hill separating them.}
\end{figure}

Once the Gibbs measure has been defined, one can define the associated Glauber dynamic --- that is, one can define a discrete-time Markov chain $(Y^N(t))_{t \in \N}$ --- which has $\mu$ as its unique invariant measure. The chain $(Y^N(t))_{t \in \N}$ on $\{0,1, \ldots, N\}$ represents the number of positive spins over time in this Glauber dynamic. For some normalizing constant $C$ and for any two contiguous integers $m,n \in \{0, 1, \ldots, N\}$, the transition probabilities are given by: $$p_{m,n} \triangleq \frac{C}{N} \sqrt{\frac{\mu(n)}{\mu(m)}} = \frac{C}{N} \sqrt{\frac{\binom{N}{n}}{\binom{N}{m}}} \, \exp\left( - \frac{1}{2T} \left(U(n) - U(m) \right) \right).$$
For non contiguous $m$ and $n$ simply set $p_{m, n} = 0$ when $|m-n|>1$ and in order for $(p_{m,n})_{0 \leq m,n \leq N}$ to be a proper stochastic matrix let $p_{n,n} = 1 - (p_{n,n-1} + p_{n,n+1})$, with the convention $ p_{N,N+1} = p_{0,-1} = 0$. An immediate calculation shows that these transition probabilities satisfy the detailed balance equation: $\mu(n) p_{n,n+1} = \mu(n+1) p_{n+1,n}$. Hence $(Y^N(t))_{t \in \N}$ has indeed the Gibbs measure $\mu$ as its unique invariant distribution. The transition probabilities can be rewritten in terms of some energy profile as follows. For any $m,n \in \{0, \ldots,N\}$ define $$\bar{a}(n) \triangleq \frac{1}{N} \left( \frac{U(n)}{T} - \log \binom{N}{n}\right)$$ and let $\Delta_{m,n}$ be the increment $\Delta_{m,n} \triangleq \bar{a}(n) - \bar{a}(m)$. Then, for $|m-n| = 1$, one can easily check that $$p_{m,n} =  \frac{C}{N} e^{- \frac{N}{2} \Delta_{m,n}}.$$
The increment $\Delta_{m,n}$ represents the amount of energy the system needs to spend for transitioning from $m$ to $n$, which is why the probability of this transition decreases exponentially fast as $\Delta_{m,n}$ increases. Finally we would like for the support of the energy profile not to be dependent on the number of particles $N$, in order to define the limit space $\mathcal{S}$ on which \ref{meta:loss} and \ref{meta:therma} can be formulated. If the system has $n$ particles with positive spin, then the average magnetization $x$ is given by $x = \frac{1}{N} \sum_{k=1}^N \sigma_k = \frac{2n - N}{N}$.
Equivalently, for a given magnetization $x$ the number $n$ of positive spins is given by $n = \frac{N}{2}(x+1)$. Therefore, if one denotes by $a_N(x)$ the energy with respect to a given magnetization $x \in \{\m 1, \m 1 + \frac{2}{N}, \ldots, 1 - \frac{2}{N},1 \} \triangleq \mathcal{S}_N$, one has $ a_N(x) = \bar{a} (n) = \bar{a} ( \tfrac{N}{2}(x+1))$.
In the limit of infinitely many particles the above function converges to some limit $a(x)$ defined on $\mathcal{S} \triangleq [\m 1,1]$ --- the Hausdorff limit of $\mathcal{S}_N$ --- with a double-well shaped graph as illustrated in Figure \ref{fig:double-well}. The function $a(\cdot)$ has a local minimum $x_1$ and a global minimum $x_2 \neq x_1$. These two minima form valleys separated by a hill, denoted by $x_0$.
This profile gives us a clear hint on how to define the metastable and stable regions. Consider the \emph{average magnetization} version of $(Y^N(t))_{t \in \N}$, denoted $(X^N(t))_{t \in \N}$ and defined for any $t \geq 0$ --- in accordance with the equations above --- as:
$$X^N(t) = \frac{2Y^N(t) - N}{N}.$$
Then $(X^N(t))_{t \in \N}$ is a Markov chain taking value in magnetization space $\mathcal{S}_N$ defined earlier. From the point of view of metastability $x_1$ is a local attractor around which $(X^N(t))_{t \in \N}$ will wander for long times, in a pseudo-stationary manner, until an unusually strong deviation leads to the crossing of the hill at $x_0$, toward stabilization around $x_2$. Hence in that context the limit metastable region is $\mathcal{M} \triangleq [\m 1,x_0)$ and $\mathcal{M}_N \triangleq \mathcal{M} \cap \mathcal{S}_N$. The limit stable region is $\mathcal{A} \triangleq (x_0,1]$. Note that the transition from $x_2$ back to $x_1$ is still possible at the mesoscopic time scale, even though it is increasingly unlikely once you've entered far enough into $\mathcal{A}$. This is why the stable region $\mathcal{A}_N$ is defined as a sufficiently deep region of $\mathcal{A}$ with respect to the volume of the system, where the possibility of crossing back the hill toward the metastable region can be safely neglected. Specifically we set $\mathcal{A}_N \triangleq [x_0 + \frac{1}{N^{1/4}}, \; 1 ] \cap \mathcal{S}_N$. Hence the tunneling time is:
$$\tau_N \triangleq \inf \left\{t \geq 0: X^N(t) > x_0 + \frac{1}{N^{1/4}}\right\}.$$
In words $\tau_N$ is the first time the system enters the absorbing region further enough for discarding the possibility of immediately crossing back the hill toward the metastable region. Throughout this example the process $(X^N(t))_{t \in \N}$ is understood to start at the point $x'_1 \in \mathcal{S}_N$ immediately to the left of $x_1$ --- this is the reference initial condition, so that the unlabeled $\P$ and $\E$ below coincide with $\P_{x'_1}$ and $\E_{x'_1}$. The macroscopic time-scale in that case is defined as $\beta_N = \E (\tau_N)$, and the mesoscopic time-scale is defined as $R_N = \beta_N/N$. Moreover, for any real-valued function $f$ and positive real number $\epsilon$, let $\Gamma^\epsilon_N(f)$ be defined as in \eqref{eq:typical}, with the integral replaced by a sum (because of the discrete nature of time in this case) and with $\mu = \delta_{x_1}$. That is, let:
$$ \Gamma^\epsilon_N(f) \triangleq \left\{\max_{\substack{t<\tau_N-R_N \\ t \equiv 0 \, \text{mod} \, R_N}}\left|\frac{1}{R_N}\sum_{s=\lfloor t \rfloor + 1}^{\lfloor t+R_N \rfloor} f\left(X^N(s)\right) - f(x_1)\right| < \epsilon \right\}.$$
From there, one can indeed rigorously prove \ref{meta:loss} and \ref{meta:therma} for the Curie--Weiss model  --- see Proposition 2.1 and Theorem 2.2 in \cite{cassandro}. That is, we indeed have asymptotic memorylessness of the hitting time $\tau_N$ on the macroscopic time-scale $\beta_N$:
$$\frac{\tau_N}{\beta_N} \overset{\mathcal{D}}{\underset{N \rightarrow \infty}{\longrightarrow}} \mathcal{E}(1),$$ where $\mathcal{E}(1)$ denotes a unit mean exponential distribution, and the superscript $\mathcal{D}$ denotes convergence in distribution under $\P$. Moreover, thermalization holds in the sense that, for any \emph{continuous} $f$ and any $\epsilon>0$ one indeed has:
$$\P \left( \Gamma^\epsilon_N(f) \right) \underset{N \rightarrow \infty}{\longrightarrow} 1.$$

Observe that while the second statement above is not completely faithful to the original formulation of Theorem 2.2 in \cite{cassandro}, it should nonetheless be clear that it is an immediate consequence from such theorem --- see the two remarks below the theorem in \cite{cassandro}. The main tool for establishing these results is a clever coupling based on the identification of some underlying renovation structure, leading to Chebyshev-type estimates. Now, for any real numbers $s>0$ and $R>0$, define the spatio-temporal average on $[s,s+R]$ as in \eqref{eq:stavg}, but again with the integral replaced by a sum:
$$ A_R^N(s, \, \cdot \,) \triangleq \frac{1}{R} \sum_{k = \lfloor s \rfloor +1}^{\lfloor s+R \rfloor} \delta_{X^N(k)} (\, \cdot \,).$$ 
Ultimately, Proposition 2.1 and Theorem 2.2 in \cite{cassandro} implies that the time-continuous and measure-valued process $(A_{R_N}^N (s \beta_N, \, \cdot \,))_{s \in \R_+}$ converges in distribution on the a suitably chosen metric space toward the measure-valued jump Markov process $( A(s, \, \cdot \, ))_{s \in \R_+}$ given by
$$A(s, \, \cdot \,) = \begin{cases}
    \delta_{x_1}(\, \cdot \,) &\text{ if } s < T,\\
    \delta_{x_2}(\, \cdot \, ) &\text{ if } s \geq T,
\end{cases}$$
where $T \sim \mathcal{E}(1)$. While the metastable distribution $\mu$ in this case is $\mu = \delta_{x_1}$, the stable equilibrium is $\nu = \delta_{x_2}$, meaning that --- starting anywhere inside the valley around $x_1$ and for big enough $N$ --- a typical trajectory of the system remains for a memory-less and exponentially long time with magnetization $x_1$, until it finally spends enough energy to reach the actual minimal magnetization $x_2$.

\subsubsection*{The Contact process}
\label{sec:contact}
The contact process is a central and somewhat canonical example of interacting particle systems ---  because it lies at the intersection of various important classes of systems due to its many regularities, much like Brownian motion is a canonical element in the general theory of stochastic processes (lying at the intersection of Markov processes, martingales, Gaussian processes and so on). Here we introduce the one-dimensional version studied in \cite{cassandro}.

\begin{figure}[ht!]
\centering 
\includegraphics[width=0.8\textwidth]{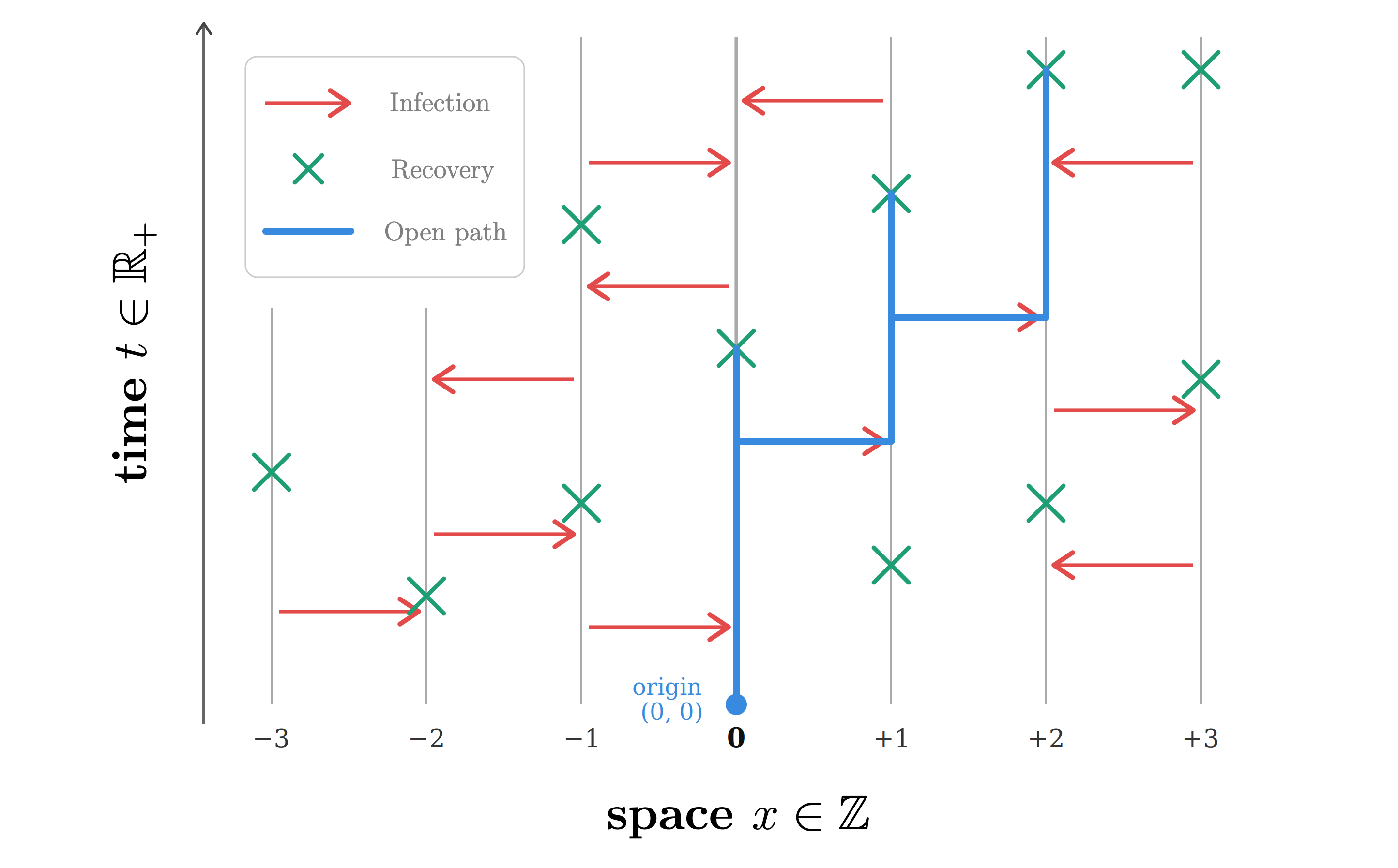}\\
\caption{\label{fig:harris} Harris graphical construction for the contact process. Here an example in which the system starts with a unique infected particle at time $0$. The red arrows are realizations of Poisson processes of rate $\lambda$ while the green crosses are realizations of a Poisson process of rate $1$. At any given time $t$ a particle is infected if and only if it can be reached by the blue path at that time.}
\end{figure}

Suppose we have a family of binary particles indexed by the set $\Z$. We call these particles ``binary'' because, much like in the Curie--Weiss model, each one of them must be in one of two different states, denoted conventionally $0$ and $1$, so that the state space of the system is $\mathcal{S} = \{0,1\}^\Z$. Typically, instead of thinking of particles with positive and negative spins, like in the previous example, one might think of these particles as individuals who might be contaminated (state $1$) or not (state $0$) by some virus. The whole dynamics is a continuous-time Markovian system modeling an epidemic, and infected individuals recover at rate $1$ while they are infected at rate  $\lambda \times \#\text{ infected neighbors}$. A common and very useful representation is the Harris percolation diagram\footnote{We refer to \cite{harris} for details.}, in which the trajectories of the system are represented on an infinite family of timelines $\Z \times \R_+$, each particle being associated with three Poisson processes, one of rate $1$ representing the recoveries (marked by crosses on the timeline $\R_+$ associated with this particle), and two of rate $\lambda$, representing the infections coming from the two neighbors on the left and on the right respectively (marked by arrows between the two timelines), as illustrated in Figure \ref{fig:harris}.
A famous result establishes that the system undergoes a phase transition with respect to $\lambda$, with a critical value $0 < \lambda_c < \infty$ such that if $\lambda < \lambda_c$ the infection dies out, whereas if $\lambda > \lambda_c$ the infection goes on for eternity. Moreover it is immediately clear that $\delta_\mathbf{0}$ --- the Dirac measure concentrated on $\mathbf{0} \triangleq 0^\Z$ --- is trivially invariant for the contact process, and the phase transition can be reformulated in terms of the invariant measures of this system using a property of the process called \emph{self-duality}\footnote{One can define a natural dual process through the graphical construction by reversing time; in the case of the contact process, this dual turns out to be the contact process itself.}: 
\begin{itemize}
    \item if $\lambda < \lambda_c$ then the only invariant measure is $\delta_\mathbf{0}$,
    \item whereas if $\lambda > \lambda_c$ there exists a non-trivial invariant measure $\mu \neq \delta_\mathbf{0}$. 
\end{itemize}

Let $(X^N(t))_{t \in \R_+}$ be the finite version of the system above, taking its values in the state space $\mathcal{S}_N \triangleq \{0,1\}^{[\m N,N]\cap \Z}$, with free boundary condition\footnote{That is, no interaction with any exterior field for the particles at the border.}. Write $\mathbf{0}_N$ and $\mathbf{1}_N$ for the all-zeros and all-ones states respectively, in the finite setting. Furthermore assume $\P = \P_{\mathbf{1}_N}$ and define:
$$\tau_N \triangleq \inf \left\{ t \geq 0: X^N(t) = \mathbf{0}_N \right\}.$$ From standard considerations on irreducibility the extinction time above is well defined, in the sense that $\tau_N < \infty$ almost surely. Once the system has reached $\mathbf{0}_N$ then no further infection can occur and the system is trapped there for eternity. Thus, the absorbing region in this case is $\mathcal{A}_N \triangleq \{\mathbf{0}_N\}$, while the metastable region is simply the complementary set $\mathcal{M}_N \triangleq \mathcal{S}_N \setminus \mathcal{A}_N$. The time-scales are defined as follows: the macroscopic time-scale $\beta_N$ is the unique solution to $\P( \tau_N > \beta_N) = e^{-1}$, and the mesoscopic time-scale is defined as $R_N = \beta_N^{9/10} N^{1/10}$. Moreover, for any real-valued function $f$ and positive real number $\epsilon$, define $\Gamma^\epsilon_N(f)$ exactly as in \eqref{eq:typical}, $\mu$ being the invariant measure of the infinite system discussed earlier. Then the two characteristic properties \ref{meta:loss} and \ref{meta:therma} can indeed be established rigorously --- see Theorem 3.2 and Theorem 3.3 in \cite{cassandro}. More precisely:
 $$\frac{\tau_N}{\beta_N} \overset{\mathcal{D}}{\underset{N \rightarrow \infty}{\longrightarrow}} \mathcal{E}(1) \text{ if } \lambda > \lambda_c.$$
Note that, even though we didn't define $\beta_N$ equal to the expectation, like in the Curie--Weiss example, a fairly immediate consequence of the result above --- when combined with the Dominated Convergence Theorem --- is that $\beta_N \sim \E( \tau_N)$ as $N \rightarrow \infty$. Concerning thermalization one indeed has that, for any cylindrical function $f$ and any $\epsilon>0$, the following holds:
$$\P \left( \Gamma^\epsilon_N(f) \right) \underset{N \rightarrow \infty}{\longrightarrow} 1 \text{ if } \lambda > \lambda_c.$$

Again, the statements above are not completely faithful to Theorem 3.2 and 3.3 in \cite{cassandro}. First of all, these results were initially proven only for $\lambda$ big enough. The refinement of the proof actually establishing the results for all $\lambda > \lambda_c$ appeared later in \textcite{schonmann}. Moreover no precise formula was given in \cite{cassandro} for the mesoscopic scale $R_N$, it was only proven that such an $R_N$ must necessarily exist. The explicit growth $R_N = \beta_N^{9/10} N^{1/10}$ was exhibited later, in \cite{schonmann} as well. One of the central tools is the Harris diagram discussed earlier. This graphical construction allows extensive use of coupling techniques, in particular between various contact processes with different initial conditions or on different spatial domains. The proofs are further facilitated by some of the nice properties of the contact process, such as \emph{monotonicity}\footnote{Roughly, monotonicity means that having more infected particles in the initial configuration implies having more infected particles in the system at any time in the future.} and self-duality.

Now, for any real numbers $s>0$ and $R>0$, define the spatio-temporal average on $[s,s+R]$ as:
$$ A_R^N(s, \cdot) \triangleq \frac{1}{R} \int_{s}^{s+R} \delta_{X^N(t)} (\, \cdot \,) \, dt.$$ 
Similarly as for the Curie--Weiss, Theorem 3.2 and Theorem 3.3 in \cite{cassandro} imply that the process $(A_{R_N}^N (s \beta_N, \, \cdot \,))_{s \in \R_+}$ converges in distribution on a suitably chosen metric space toward the measure-valued jump Markov process $( A(s, \, \cdot \, ))_{s \in \R_+}$ defined by
$$A(s, \, \cdot \,) = \begin{cases}
    \mu &\text{ if } s < T,\\
    \delta_{\mathbf{0}} &\text{ if } s \geq T,
\end{cases}$$
with $T \sim \mathcal{E}(1)$. In that case, the metastable distribution is the non trivial invariant distribution $\mu$ of the infinite contact process, and the stable equilibrium is $\nu = \delta_\mathbf{0}$. This means that a typical trajectory of the finite system $(X^N(t))_{t \in \R_+}$ will initially behave in a pseudo-stationary fashion approximately described by the restriction of $\mu$ to $\mathcal{S}_N$, until an unusually large number of consecutive recoveries abruptly leads the system to absorption in $\mathbf{0}_N$.

\section{The Galves--Löcherbach models} \label{sec:GLModels}

The basic biological picture goes as follows. The elementary component is the neuron. Besides the cell body (the soma), a neuron is composed of an axon --- its ``output wire'' --- and many dendrites --- its ``input branches''. Neurons receive signals from other neurons through the synapses --- the junctions between a neuron's axon and another neuron's dendrites. These incoming signals activate the ion channels on the membrane of the neuron, modifying the \emph{membrane potential} of this cell --- a crucial quantity defined as the difference in electrical potential between the interior and the exterior of the cell. As the membrane potential rises, the probability increases that the neuron emits an action potential --- or \emph{spike} --- a brief, stereotyped electrical pulse that travels along the axon. When it reaches the synapses, the spike triggers the release of neurotransmitters into the synaptic cleft, which in turn influences the membrane potential of all neurons downstream. The synaptic connection can be either excitatory (pushing the target neuron closer to firing) or inhibitory (pushing it further away). After firing, the neuron's potential resets to a resting value, so that the future evolution of its membrane potential depends only on the activity received since this last spike. Moreover, membrane potentials do not stay fixed between incoming spikes --- they decay spontaneously toward a resting value due to ion channel dynamics, a phenomenon known as \emph{leak currents}. 

The overall activity of a neural circuit thus emerges from this constant interplay between excitation, inhibition, integration, leakage, and reset. Spikes are widely regarded as the fundamental unit of neural information --- roughly similar to what bits are to computers. We shall mention that this spiking mechanism --- transmitting digital-like information through discrete events, and mediated by neurotransmitter release --- corresponds to communication through chemical synapses, which are by far the most numerous and functionally dominant type of synapse in the brain. There also exist electrical synapses, allowing a continuous and bidirectional current flow between neurons, though these play a more specialized role in neural circuits.

The mathematical modeling of neurons and neural networks now has a long history, since the introduction of the classical integrate-and-fire model by Louis Lapicque at the beginning of the 20th century. However, more recent findings indicate that empirical evidence is at odds with the deterministic views of the original model and with most of the more refined models it has inspired (such as the famous Hodgkin--Huxley model), all based on the same mathematical tool --- i.e.\@ differential equations. Indeed stochasticity intervenes at many levels in the dynamics of neural networks: ion channels open and close randomly, synaptic transmission involves random release of neurotransmitter quanta, and the same stimulus does not reliably produce exactly the same spike pattern twice. Hence biological neural networks are essentially big systems of random interacting components. These components are individually relatively well understood, but still poorly understood at the macroscopic level of the whole system, mostly because of the astronomical number of components (of the order of $10^{11}$) and because of the complexity of their interactions. This is typically the kind of situation statistical physics is good at handling, which suggests elaborating a model amenable to its methods.

\subsection{General framework}

These considerations led to the elaboration of a family of stochastic models --- generically called GL models --- aiming at giving a mathematically tractable description of neural networks incorporating all of the biological features mentioned above, by describing each neuron as an individual stochastic unit rather than resorting to the intractable precision of detailed biophysical models of Hodgkin--Huxley type. The original GL model, introduced in \cite{glmodel}, is a discrete-time stochastic process, later extended to continuous-time in \textcite{gl4}. This continuous-time version incorporated features such as the emission and transmission of spikes through chemical synapses, and continuous drift toward an average value through electrical synapses, but did not yet include leak currents. The model was therefore further generalized in \textcite{duarteost} to incorporate this phenomenon. Other continuous-time versions, sharing the framework of \cite{duarteost} but adopting a different leakage mechanism, were subsequently introduced in \textcite{ferrari} and \textcite{amarcos}. It shall be noted that the fundamental ideas behind this model makes it relatively versatile and that it can actually apply outside the realm of biological neural networks --- similar ideas have been considered for modeling opinion dynamics in \cite{galveslaxa,laxa} for example, in a framework very closely related to the one introduced below. Since metastability is more naturally studied in a continuous-time setting --- especially from the point of view of the pathwise approach introduced earlier --- we will focus below on these continuous-time versions of the model. A self-contained, textbook-level development of this whole program can be found in the recent monograph \cite{gb}.

We consider a countable set of neurons $I \subset \Z$, which can be finite or infinite. For each pair of neurons $i,j \in I$, the \term{synaptic weight} of the synapse from $i$ to $j$ is denoted $w_{i \rightarrow j} \in \mathbb{R}$. These weights account for the strength of the connection, but also for its nature --- either excitatory or inhibitory, depending on whether $w_{i \rightarrow j} > 0$ or $w_{i \rightarrow j} < 0$ respectively. We will generically assume that $w_{i \rightarrow i}=0$ for all $i \in I$ --- that is, there is no self excitation nor self inhibition --- and we will furthermore assume the following uniform summability condition to avoid pathological cases\footnote{This condition can actually be weakened, since by considering only absolute values we neglect the possibility of a balance between excitation and inhibition, see \cite{papageorgiou}.}
$$
\sup_{i \in I}\sum_{j \in I}|w_{j \rightarrow i}| < +\infty.
$$

From the high level description viewpoint, the continuous-time GL model is a system of interacting point processes. The spike train of each neuron $i \in I$ is described by its own point process $Z_i$ on the time line $\R_+$. For any $s<t \in \R_+$, the number of spikes emitted by neuron $i$ on the time interval $[s,t]$ is denoted $Z_i([s,t])$, and $Z_i(t)$ is an abbreviation for $Z_i([0,t])$. Unlike traditional Poisson point processes the intensity $(\lambda_i(t))_{t \in \R_+}$ for any given neuron $i$ is itself a stochastic process, and unlike Cox processes the stochasticity of $(\lambda_i(t))_{t \in \R_+}$ is endogenous to the system, meaning that the value of $\lambda_i(t)$ at any given time $t$ is directly influenced by the values of $Z_j(s)$, for $s \leq t$ and $j \in I$. Formally,  let $\mathcal{F}_t \triangleq \sigma(Z_i(s): 0 \leq s \leq t, i \in I)$ be the filtration of the whole system up to time $t$. Then, for all $i \in I$, we have:
\begin{enumerate}
\item $\P \left(Z_i([t,t+dt])=1 \, \middle| \, \mathcal{F}_t \right) = \lambda_i(t)dt + o(dt)$,
\item and $\P \left(Z_i([t,t+dt]) \geq 2 \, \middle| \, \mathcal{F}_t\right) = o(dt)$.
\end{enumerate}
In contrast to linear Hawkes processes, where the intensity is an affine function of past events, the \term{spiking rate} $\lambda_i(t)$ of neuron $i$ depends on its \term{membrane potential} $U_i(t)$ --- to be defined in subsequent sections --- as follows:
$$
\lambda_i(t) = \phi\left(U_i(t)\right),
$$
where $\phi$ is a non-decreasing real function called the \term{activation function}. Quite crucially, unlike nonlinear Hawkes processes --- a popular model for biological neural networks, which essentially fits the description up to now --- the dependence of the intensity $\lambda_i(t)$ on the past activity of the system is reset at each firing event of neuron $i$, making the intensity depend only on the activity received since its last spike. This means that each neuron has \emph{variable length memory}, and this property is what fundamentally distinguishes the GL model from Hawkes-type processes. This specificity of the model is introduced at the level of the membrane potentials $U_i(t)$, where the interplay between excitation, inhibition, leak currents, and resets particular to neural networks takes place. Below we distinguish between three different closely related possibilities for the definition of this membrane potential, which differ by the way they model leak currents. Since the versions of the GL model reviewed in Section \ref{sec:metastability_GLmodel} are all purely excitatory --- i.e.\@ they satisfy $w_{i \rightarrow j} > 0$ for all pairs $i \neq j$  --- we'll assume below that $U_i(t) \in \R_+$ for all $i$ and $t$, even though it doesn't have to be the case in all generality.

\subsection{Continuous leakage}

The first and most straightforward option for incorporating leak currents is to model them as a continuous and deterministic decay of the influence of past spikes on the membrane potential, represented by a decreasing and measurable \term{leak function} $g: \mathbb{R}_+ \to \mathbb{R}_+$. For each neuron $i \in I$ and each time $t \in \mathbb{R}$, define the last spiking time of neuron $i$ before time $t$:
$$ L_i(t) \triangleq \sup\{s < t : \Delta Z_i(s) > 0\},$$
where $\Delta Z_i(s) \triangleq Z_i(s) - Z_i(s^-)$. The membrane potential of neuron $i$ at time $t$ is then defined as:
$$U_i(t) \triangleq \sum_{j \in I} w_{j \rightarrow i} \int_{[L_i(t),\, t)} g(t - s)\, dZ_j(s).$$
At each spiking time of neuron $i$, its membrane potential is reset to $0$, while the membrane potential of any other neuron $j \in I \setminus \{i\}$ is shifted by $w_{i \rightarrow j}$. This reset causes the memory of neuron $i$ to have variable length: its membrane potential depends only on the activity of the other neurons since its last spiking time. Furthermore, between any two consecutive spikes of the network, the membrane potential of each neuron evolves deterministically according to $g$. The interpretation of $g$ is the following: if a given neuron $j$ emitted a spike $t$ time units ago, its current influence on the membrane potential of any other given neuron $i$ amounts to $g(t)$ times its original value.

Besides the trivial case in which there is no leakage ($g \equiv 1$) there is exactly one case in which the system of membrane potentials --- that is, the multidimensional process  $(U(t))_{t \in \R_+} = (U_i(t))_{i \in I, \,t \in \R_+}$ --- is actually Markovian. This is the case of exponential leakage: $g(t)=e^{-\alpha t}$ for some $\alpha >0$. Then $(U(t))_{t \in \R_+}$ is a piecewise deterministic Markov process\footnote{See \cite{davis} for a precise definition of this class of Markov processes.} with infinitesimal generator
\begin{equation}\label{eq:GLmodelexp_leakage}
\mathcal{L}f(u)=\sum_{i\in I}\phi(u_i)[f(\pi_i(u))-f(u)] - \alpha \sum_{i\in I} \frac{\partial f}{\partial u_i}(u)u_i,
\end{equation}
for any $u \in \mathbb{R}_+^I$ and any smooth test function $f:\mathbb{R}_+^I \to \mathbb{R}$, the maps $\pi_i$ being defined by:
\begin{equation} \label{eq:spikingmaps}
\left(\pi_i(u) \right)_j \triangleq
\begin{cases}
u_j+w_{i \rightarrow j}, &\text{ if } j\neq i, \\
0 &\text{ if } j= i.
\end{cases}
\end{equation}

In words, at each spiking time of a given neuron $i \in I$, its own membrane potential resets to $0$ and the membrane potential of any other neuron $j \in I \setminus\{i\}$ increases by an amount of $w_{i \rightarrow j}$ units. Moreover, in between two consecutive spikes of the system, the membrane potential of each neuron $i\in I$ approaches $0$ deterministically at exponential speed with rate $\alpha$. This version of the model is exactly the one introduced in \cite{duarteost}, once electrical synapses have been neglected --- which is the typical setting in which metastability is expected. Assuming exponential leak is not only mathematically convenient --- as it places us in a Markovian setting --- but it is also quite realistic biologically, being the standard format for leak currents adopted in traditional models in computational neuroscience, such as leaky integrate-and-fire models inspired by Lapicque's work or the Hodgkin--Huxley model.

\subsection{Total leakage}
The exponential leakage version of the model introduced above has the advantage of being still quite realistic on a biophysical ground. But biophysical realism in the dynamical properties of the system is not free: it usually comes at the expense of additional mean-field type assumptions on the topology of the network in order to reach mathematical tractability. Moreover one could argue that --- as long as the effect of the interplay between spikes and leakage on the \emph{timing of the spikes} is included in the model one way or another --- the exact trajectory of the membrane potential in between two spikes isn't so much relevant to us. For these reasons, simpler alternatives have been considered for introducing leakage into the model.

One of these alternatives, is the use of \term{total leakage} times, modeling roughly via punctual events the time it takes for a neuron to reach membrane potential $0$. Formally, suppose we associate to each neuron a Poisson point process $(Z_i^\dagger(t))_{t \in \R_+}$ of some rate $\gamma > 0$. Redefine the memory length of neuron $i$ by:
$$L_i(t) \triangleq \sup\{s < t : \Delta Z_i(s) > 0 \text{ or } \Delta Z_i^\dagger(s) > 0\},$$
and define the membrane potential as usual but without continuous leaks --- i.e.\@ $g(\cdot) \equiv 1$:
$$U_i(t) = \sum_{j \in I} w_{j \rightarrow i} \int_{[L_i(t),\, t)}\, dZ_j(s).$$
The system of the membrane potentials $(U(t))_{t \in \R_+}$ is now a Markov jump process with infinitesimal generator:
\begin{equation}\label{eq:GLmodeltotalleakage}
\mathcal{L}f(u)=\sum_{i\in I}\phi(u_i)[f(\pi_i(u))-f(u)]+ \gamma \sum_{i\in I}[f(\pi_i^{\dagger}(u))-f(u)],
\end{equation}
for any $u \in \mathbb{R}_+^I$ and any smooth test function $f:\mathbb{R}_+^I \to \mathbb{R}$. Here the spiking maps $\pi_i$ are still defined as in \eqref{eq:spikingmaps} and the leakage maps $\pi_i^\dagger$, are defined by:
\begin{equation}
\left(\pi_i^\dagger(u) \right)_j \triangleq u_j \mathbbm{1}_{i \neq j}
\end{equation}

In words, in that alternative modeling strategy, the time at which a neuron loses all potential --- all other things being equal --- is modeled directly by a random variable with exponential distribution of rate $\gamma$ --- instead of resulting secondarily from a deterministic decay with exponential speed of rate $\alpha$. The mean inter-leak time $1/\gamma$ might thus be seen as the average time it takes for a neuron starting from a ``typical'' value of membrane potential and decaying deterministically at a given rate $\alpha$, to reach a potential near $0$. The spiking mechanism on the other hand remains exactly as in the preceding version of the model. This version of the model was introduced by \cite{ferrari} in the specific case of a binary activation function $\phi(u) = \mathbbm{1}_{u > 0}$. As we will see later, one of the major advantages of modeling leak currents via total leakage times is that the analysis is significantly simplified by allowing a thorough study of the IPS keeping track of the spiking rates.

\subsection{Unitary leakage}
A third leakage mechanism still models leaks via punctual events, but instead of an 
immediate reset to the resting value $0$, each such event provokes a discrete unit 
decrement of the membrane potential. A natural further generalization, motivated by 
analogy with the continuous leakage setting, is to allow the rate of these leakage 
events to depend on the membrane potential itself. Formally, let $\alpha : \mathbb{R}_+ \to \mathbb{R}_+$ be the rate function 
of the leaks. Then $(U(t))_{t \in \R_+}$ is a Markov jump process with 
infinitesimal generator
\begin{equation}
\mathcal{L}f(u)=\sum_{i\in I}\phi(u_i)[f(\pi_i(u))-f(u)]+ \sum_{i\in I} 
\alpha(u_i)[f(\pi_i^\dagger(u))-f(u)],
\end{equation}
for any $u\in \mathbb{R}_+^I$, where the spiking maps $\pi_i$ are still defined as 
in~\eqref{eq:spikingmaps} and the unit leak maps $\pi_i^\dagger$ are defined by
\begin{equation}
\left(\pi_i^\dagger(u)\right)_j \triangleq \left( u_j - \mathbbm{1}_{i = j} \right) \vee 0.
\end{equation}
This version of the model was introduced by \cite{amarcos} in the natural case 
$\alpha(u) = \gamma u$ for some $\gamma > 0$, which makes the leakage rate proportional 
to the current membrane potential, and again with binary activation function $\phi(u) = \mathbbm{1}_{u > 0}$.

\section{Metastability in Galves--Löcherbach models} \label{sec:metastability_GLmodel}

The results we survey differ along several essentially independent axes: the leakage mechanism (continuous, total, or unitary), the activation function $\phi$ (binary, saturated linear, or exponential), the interaction graph (complete/mean-field, lattice, or otherwise), and the precise approach which is considered --- pathwise or otherwise --- and in the pathwise case, which of the two properties \ref{meta:loss} and \ref{meta:therma} is actually established. Table~\ref{tab:summary} summarizes this landscape. 

We mainly focus on the works in which metastability is established rigorously in the pathwise sense, giving an in-depth presentation along the three first subsections, and sketching the proofs when possible. Other works are discussed separately and more briefly in the fourth subsection, either because they do not completely fit the GL framework we defined in the previous section, or because they are of numerical nature or less clearly related to our pathwise perspective.

As we sketch the proofs, we try to shed light on the structural similarities, in order to make the general method explicit. We will especially focus on the proof of \ref{meta:loss}, since few works among the ones reviewed here have explicitly worked out \ref{meta:therma} using the precise formalism of the pathwise approach introduced earlier. The starting point usually consists in establishing the relevant upper-bound on the \emph{factorization defect} of the extinction time, in order to prove that it vanishes as the volume of the system diverges. As we will see a common key point in all these proofs is the partition of the metastable region into two subsets. The subset $\mathcal{W}_N \subset \mathcal{M}_N$ represents the \term{typical states} --- i.e.\@ the ones in the inner region of $\mathcal{M}_N$, where any typical trajectory of the system is expected to spend most of its time --- and where one can generally identify some renovation pattern. The complementary set $\mathcal{M}_N \setminus \mathcal{W}_N$ then represents the \term{atypical states} --- the ones at the margin of the metastable region, which are vanishingly unlikely in typical trajectories of the system. The bound obtained on the factorization defect has generally two terms, related to the typical and atypical states respectively. The techniques underlying the proofs that these two terms vanish are another common key point to all the works described in this section. Indeed, each case involves at some point a kind of \emph{synchronous coupling}, which is used to prove that two copies starting from two different points in $\mathcal{W}_N$ are expected to collide quickly, making the system effectively lose its memory. Moreover, in each case one has also to consider a \emph{proxy} for the system at some point --- a modification of the initial system or a mean-field limit for example --- having an invariant measure which can be either explicitly computed or proven to concentrate on specific regions of the state space. We will generally denote by $\mu$ this invariant measure, which approximately describes the distribution of the original system in the metastable region before tunneling. In some cases it is literally the same mathematical object as the invariant measure $\mu$ introduced in Section \ref{subsec:pathwisedef}, while in other cases it merely plays an analogous role with some mathematical subtlety --- in particular we might have $\mu = \mu_N$ dependent on the volume and defined on $\mathcal{S}_N$ instead of $\mathcal{S}$.

\textbf{Notation.} From now on, since metastability occurs in finite systems, the stochastic processes involved will frequently be written with a superscript $N$ referring to the volume --- i.e.\@ the number of neurons in the system --- denoting for example $(U^N(t))_{t \in \R_+}$ instead of $(U(t))_{t \in \R_+}$, $(Z_i^N(t))_{t \in \R_+}$ instead of $(Z_i(t))_{t \in \R_+}$ and so on. The dependence on the \emph{initial condition} is, on the other hand, carried by the probability rather than by the process: we write $\P_x$ for the law of a process started from the state $x$ --- or equivalently $\P_\nu$ when the initial state is drawn from a distribution $\nu$ --- while an unlabeled $\P$ refers to a reference initial condition fixed by the context --- typically the \emph{maximal configuration}, in whatever sense is relevant there. In particular a hitting time such as $\tau_N$ denotes a single random variable, the same measurable functional of the trajectory regardless of where the process starts; its dependence on the starting point being carried through $\P_x$. We consistently apply similar conventions for the expectation. The only exception occurs when two coupled copies of the system must appear under the same probability. In that last case, the hitting time of the copy started from a configuration $x$ may instead be labeled by a superscript --- e.g.\@ $\tau_N^x$ --- the unlabeled $\tau_N$ still corresponding to the reference initial condition. 

\begin{table}[ht!]
\centering
\footnotesize
\renewcommand{\arraystretch}{1.25}
\begin{tabular}{llllll}
\hline
Reference & Leakage & Activation $\phi$ & Interaction & Parameter & Result \\
\hline
\cite{mandre1}     & total            & $\mathbbm{1}_{u>0}$        & $\Z$, n.n.   & $N\to\infty$      & \ref{meta:loss} \\
\cite{mandre3}     & total            & $\mathbbm{1}_{u>0}$        & $\Z$, n.n.   & $N\to\infty$      & \ref{meta:therma} \\
\cite{mandre2}     & total            & $\mathbbm{1}_{u>0}$        & complete     & $N\to\infty$      & \ref{meta:loss} \\
\cite{romaro}      & total            & $\mathbbm{1}_{u>0}$        & $\Z^d$, n.n. & $N\to\infty$      & numerical \\
\cite{laxa}        & total / unitary  & $e^{u}\mathbbm{1}_{u>0}$   & complete     & $N\to\infty$      & \ref{meta:loss} \\
\cite{evamonm}     & exponential      & $\lambda u \wedge \phi^*$ & mean-field   & $N\to\infty$      & \ref{meta:loss} \\
\cite{galveslaxa}  & none (reset)     & exponential               & complete     & $\gamma\to\infty$ & \ref{meta:loss} \\
\cite{pouzatandre} & none             & $\mathbbm{1}_{u\geq\theta}$& complete     & $N\to\infty$      & QSD / numerical \\
\cite{taille}      & exponential      & $h_ie^{\kappa_i u}$       & two groups   & $R\to\infty$      & numerical (RMF) \\
\cite{ludmiguel}   & geometric        & $u\wedge 1$               & mean-field   & $N\to\infty$      & numerical / bounds \\
\hline
\end{tabular}
\caption{\label{tab:summary} Overview of the metastability results for Galves--Löcherbach-type models surveyed in Section~\ref{sec:metastability_GLmodel}. Here \ref{meta:loss} refers to the asymptotic loss of memory and \ref{meta:therma} to thermalization; ``n.n.'' stands for nearest-neighbor and ``RMF'' for replica-mean-field. The parameter column records the limit under which the result is stated. Entries marked ``numerical'' are supported by simulations rather than by a rigorous proof.}
\end{table}

\subsection{Binary activation function and total leakage} \label{subsec:metabinaryactivationfunction}

The version of the GL model with total leakage introduced in \cite{ferrari} was initially studied in the particular case of a one-dimensional lattice with nearest neighbors interaction, under the perspective of phase transition. It was then studied under the perspective of metastability in \cite{mandre1,mandre2,mandre3,romaro} both in the original lattice setting and in a complete-interaction setting. The main underlying idea of these works is to further simplify the model by taking a hard-threshold for the activation function: $\phi(u) = \mathbbm{1}_{u>0}$. One can then focus on the auxiliary system defined by the spiking rates over time $$(\lambda^N(t))_{t \in \R_+} = (\lambda^N_1(t), \lambda^N_2(t), \ldots, \lambda^N_N(t))_{t \in \R_+},$$ where $\lambda^N_i(t) = \mathbbm{1}_{U^N_i(t) >0} \in \{0,1\}$. Then $(\lambda^N(t))_{t \in \R_+}$ is a Markovian process in $\{0,1\}^N$ and can therefore be seen as an interacting particle system and studied as such. Below we go against the chronological order by presenting the complete-interaction version first, because the lack of spatial structure makes it actually less complicated than the lattice version.

\subsubsection*{Complete interaction}

The same version of the model has been studied under complete interaction --- i.e.\@ $w_{i \rightarrow j} = 1$ for all neurons $i,j \in I$ --- in \textcite{mandre2}, establishing in particular \ref{meta:loss}. The complete interaction somewhat simplifies the proof. One indeed only needs to consider the count of active neurons over time in that case, instead of the intricate configurations of $\{0,1\}^N$, reducing the problem to the study of a relatively simple \term{aggregated process} on $\mathcal{S}_N \triangleq \{0,1, \ldots, N\} \subset \N \triangleq \mathcal{S}$. More precisely, if $S^N(t) \triangleq \sum_{k=1}^N \lambda^N_k(t)$ then, thanks to the absence of spatial structure, $(S^N(t))_{t \in \R_+}$ is a Markov jump process on $\mathcal{S}_N$ with intensity matrix $Q = ( q_{i,j})_{0 \leq i,j \leq N}$ given, for any $i \neq j $, by:
$$q_{i, j} = \begin{cases}
    \gamma i &\text{ if } i \neq 0 \text{, } i \neq N \text{ and } j=i-1,\\
    i &\text{ if } i \neq 0 \text{, } i \neq N\text{ and } j=N-1,\\
    (1+\gamma) N &\text{ if } i = N \text{ and } j=N-1,\\
    0 &\text{ otherwise.}
\end{cases}$$
These transition rates deserve a word of explanation. For the second case, under complete interaction a single spike resets the firing neuron and, simultaneously, pushes every currently sub-threshold neuron above threshold, so that the post-spike count is always $N-1$, irrespective of the number $i$ of active neurons just before the spike. Since each of the $i$ active neurons fires at rate $1$, this transition occurs at rate $i$. The leak events, occurring at rate $\gamma$ per active neuron, account for the downward transitions $i \to i-1$ at rate $\gamma i$ --- and, from the full configuration $i=N$, both a spike and a leak lead to $N-1$, whence the rate $(1+\gamma)N$.

The extinction time is $$\tau_N \triangleq \inf \left\{t  \geq 0: S^N(t) = 0 \right\}.$$ Hence the absorbing region is $\mathcal{A}_N \triangleq \{0\}$ and the metastable region is simply $\mathcal{M}_N \triangleq \mathcal{S}_N \setminus \mathcal{A}_N$. Following the convention stated earlier the unlabeled $\P$ denotes the law of $(S^N(t))_{t \in \R_+}$ when it starts from the full initial condition $S^N(0) = N$. Property \ref{meta:loss} is established as an immediate consequence of: 
\begin{equation*}
\P \left( \frac{\tau_N}{\beta_N} > s + t \right) \underset{N \rightarrow \infty}{\longrightarrow} \P \left(\frac{\tau_N}{\beta_N} > s\right) \, \P \left( \frac{\tau_N}{\beta_N} > t \right).
\end{equation*}
By elementary computations one can establish the following upper-bound on the factorization defect:
\begin{multline*}
\left| \P \left( \frac{\tau_N}{\beta_N} > s + t \right) - \P \left( \frac{\tau_N}{\beta_N} > s \right) \P\left(  \frac{\tau_N}{\beta_N} > t \right) \right| \\
\leq \sum_{k \in \mathcal{W}_N} \left| \P_k \left( \frac{\tau_N}{\beta_N} > t \right) - \P \left( \frac{\tau_N}{\beta_N} > t \right) \right| + \sum_{k \notin \mathcal{W}_N} \P (S^N(\beta_N s) = k),
\end{multline*}
where $\mathcal{W}_N \triangleq [\frac{N}{2}, N] \cap \N$. This set designates the higher end of the metastable region --- on which the typical trajectories of the system concentrate most of their life-time --- and one shall use some kind of renovation argument to show that, if $k \in \mathcal{W}_N$, then $\tau_N$ has nearly the same law under $\P_k$ and under the full start, in a sense that is made precise below, implying that the corresponding sum in the upper-bound vanishes as $N$ diverges. The lower end of the metastable region $\mathcal{M}_N \setminus \mathcal{W}_N$ contains only atypical states, rarely visited by the typical trajectories, and one needs to prove that they are unlikely enough for the corresponding sum in the upper-bound to vanish as well when $N$ diverges.

To show that the sum over $\mathcal{W}_N$ goes to $0$ as $N$ diverges one might consider a joint construction on the same probability space --- a coupling --- of two copies of the aggregated process, one started from $k$ active neurons and the other from the full configuration $N$, on which both reach $N-1$ together at the first spike of the $k$-started copy and merge into a single trajectory after that point --- this is the synchronicity allowing the renovation argument we mentioned above   --- unless $k$ leaks occur one after the other before any spike, which has probability $(\frac{\gamma}{1+\gamma})^k$. Since we're summing over the $\frac{N}{2} + 1$ elements of $\mathcal{W}_N$ this last eventuality occurs with probability less than $(\frac{N}{2} + 1)(\frac{\gamma}{1+\gamma})^{N/2}$, and hence the sum goes quickly to $0$ as $N$ diverges.

For the sum over $\mathcal{M}_N \setminus \mathcal{W}_N$ observe that, while $(S^N(t))_{t \in \R_+}$ has only one invariant measure, concentrated on $0$, one might consider an artificial modification of this process where the transition to $0$ has been removed, and couple it with the original process in such a way that they coincide for as long as possible --- by making all transition intensities equal, except the transition from $1$ to $0$, which becomes null, and the transition from $1$ to $N-1$, which receives the value $1 + \gamma$. This modified version becomes a proxy for the system itself, which obviously dominates the original version, and it has a unique non-trivial invariant measure $\mu$ which can be computed explicitly. The state $N$ is definitively left after the first spike/leak, so that $\mu_N=0$, and for any $1 \leq k \leq N-1$, one has:
$$\mu_k \propto \frac{1}{k}\left(\frac{1 + \gamma}{\gamma} \right)^{k-1}.$$
Hence the mass under stationary regime concentrates (almost) exponentially fast on the typical states --- since there is a factor of at least $\frac{1}{N}(\frac{1+\gamma}{\gamma})^{N/2}$ between the typical and atypical states. In other words --- under stationary regime, and modulo some negligible linearly decreasing factor --- the mass on the atypical states vanishes exponentially fast as $N$ diverges. Now notice that since $\max_{1 \leq k \leq N} \mu_k = \mu_{N-1}$ one necessarily has $\mu_{N-1} \geq 1/N$ and therefore $$\P_{N-1} \left(S^N(t) = k\right) \leq \frac{1}{\mu_{N-1}} \P_\mu \left( S^N(t) = k \right) \leq N \, \P_\mu \left( S^N(t) = k \right).$$
Hence the mass on the lower end of the state space vanishes exponentially fast as well for the system starting from $N-1$, even though it is not in stationary regime, since in the worst case this cannot add more than another irrelevant linear factor. This proves the result for the system starting at $N-1$, and the result for the system starting at $N$ follows immediately, since after a vanishingly small delay, $(S^N(t))_{t \in \R_+}$ always jumps to $N-1$ anyway.

\subsubsection*{Spatially structured interaction}
Now we turn to the lattice model. The neurons are indexed by $\Z$ or a finite window of contiguous elements of $\Z$ and the system has nearest-neighbors interaction, meaning that $w_{i \rightarrow j} = 1$ if $|i-j|=1$, and $w_{i \rightarrow j} = 0$ otherwise. This model was proven to exhibit a phase transition with respect to the leakage parameter in \cite{ferrari}, in the infinite setting. There exists a critical value $0 < \gamma_c < \infty$ such that: if $\gamma > \gamma_c$ the only invariant measure for the infinite system is the Dirac measure concentrated on the null state $0^\Z$, whereas if $\gamma < \gamma_c$ there exists another \emph{non-trivial} and \emph{spatially ergodic} invariant measure $\mu$ for the infinite system.

\begin{figure}[ht!]
\centering 
\includegraphics[width=0.8\textwidth]{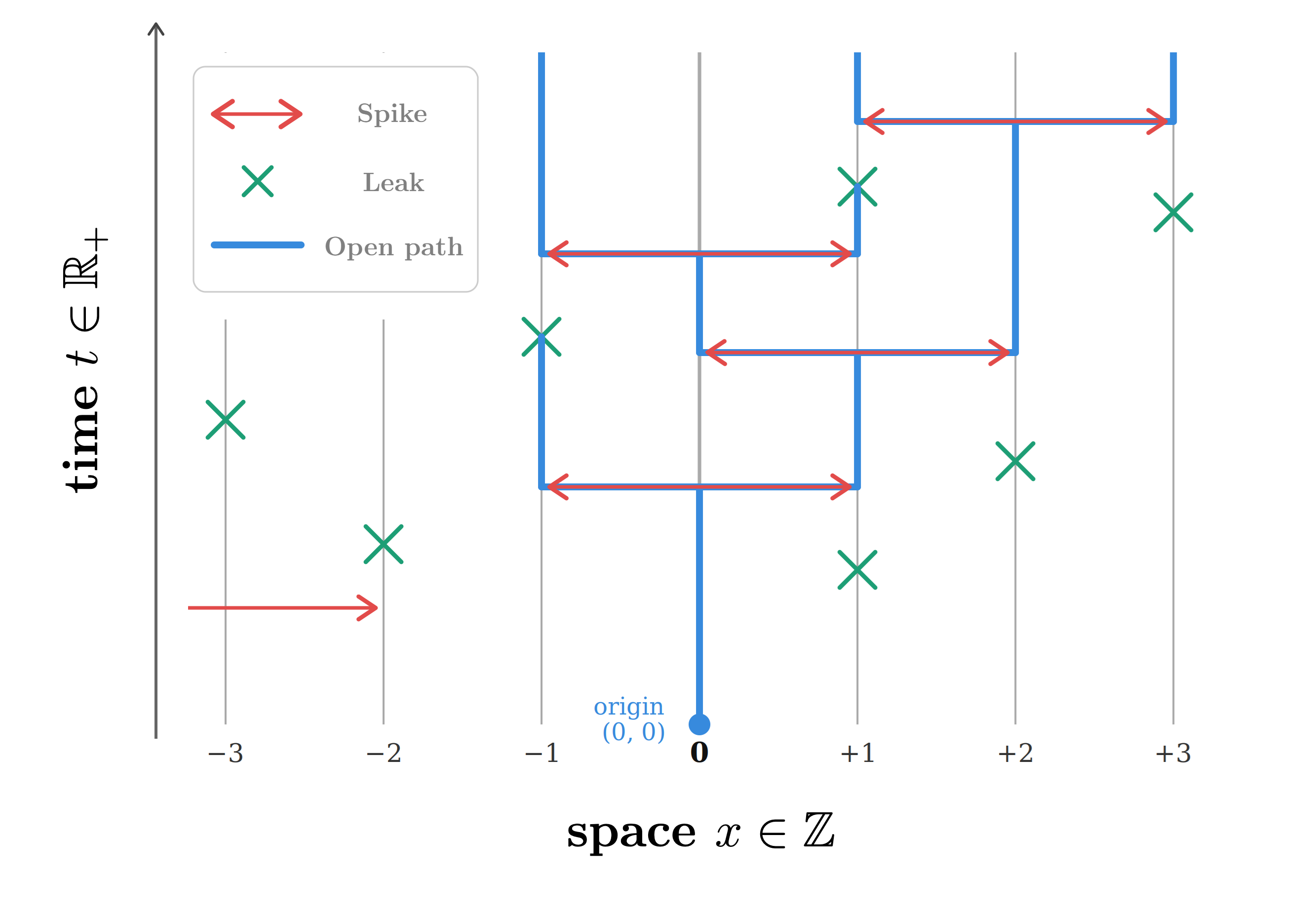}\\
\caption{\label{fig:harris_spiking} Harris graphical construction for the binary spiking neuronal network. Here an example in which the system starts with a single active neuron at the origin at time $0$. The red double arrows are realizations of independent Poisson processes of rate $1$, representing spikes --- which affect both neighbors simultaneously --- while the green crosses are realizations of independent Poisson processes of rate $\gamma$, representing total leaks. At any given time $t$ a neuron is active if and only if it can be reached by the blue path at that time.}
\end{figure}

The model is reminiscent of the one-dimensional version of the contact process discussed in Section \ref{sec:contact}. In particular, both processes share common desirable properties, such as monotonicity, nearest-neighbors interaction and translation invariance. As such, the metastable properties of $(\lambda^N(t))_{t \in \R_+}$ can be studied using similar techniques to those used in \cite{cassandro} in the study of the metastable properties of the contact process. In particular, the system can be constructed via Harris graphical representation in a similar way to the contact process itself --- see Figure \ref{fig:harris_spiking} --- which gives us access to similar coupling techniques. 

A central difference between the two processes nonetheless is that, in the neuronal model, transitions occur at multiple coordinates at the same time, while they occur only one coordinate at a time in the contact process. Indeed whenever a spike occur in $(\lambda^N(t))_{t \in \R_+}$, both neighbors are ``infected'' (or activated) simultaneously, while only one neighbor is infected at a time in the contact process. Moreover, when emitting a spike, the neuron deactivates itself (its value goes from $1$ to $0$) while the particles of the contact process remain infected when they infect a neighbor. This moves $(\lambda^N(t))_{t \in \R_+}$ out of a class of interacting particle systems sometimes generically referred to as \emph{spin systems} (see Chapter 3 of \cite{liggett}). In particular, unlike the contact process, $(\lambda^N(t))_{t \in \R_+}$ is not self-dual. One can nevertheless still adapt the techniques used in \cite{cassandro} at the price of some additional technical difficulties --- addressed in \cite{mandre1, mandre2}.

As for the contact process setting the general state space is $\{0,1\}^\Z$, or equivalently $\mathcal{S} \triangleq \mathcal{P}(\Z)$\footnote{In this section we deliberately abuse notation by identifying \emph{class of sets} and \emph{sequences of 0 and 1}, which is a frequently used convention in IPS literature.}--- the class of all subsets of $\Z$ --- and $(\lambda^N(t))_{t \in \R_+}$ will denote the process living in $\mathcal{S}_N \triangleq \mathcal{P}(\Z \cap [\m N,N])$. The full start $\lambda^N(0) = \Z \cap [\m N, N]$ is the reference initial condition, corresponding to the unlabeled $\P$ (and to the unlabeled extinction time $\tau_N$ below). The stable region will again be $\mathcal{A}_N \triangleq \{\emptyset\}$ and the rest of the state space is the metastable region $\mathcal{M}_N \triangleq \mathcal{S}_N \setminus \mathcal{A}_N$. Define the extinction time of the system as:
$$\tau_N \triangleq \inf \left\{t  \geq 0: \lambda^N(t) = \emptyset \right\}.$$
The memorylessness property \ref{meta:loss} was established in \textcite{mandre1} for a macroscopic time $\beta_N$ defined indifferently by $\P( \tau_N > \beta_N) = e^{-1}$ or $\beta_N = \E( \tau_N )$. For $\gamma < \gamma_c$ we indeed have: $$ \frac{\tau_N}{\beta_N} \overset{\mathcal{D}}{\underset{N \rightarrow \infty}{\longrightarrow}} \mathcal{E} (1).$$
As usual, this result is established as consequence of
\begin{equation*}
\P \left( \frac{\tau_N}{\beta_N} > s + t \right) \underset{N \rightarrow \infty}{\longrightarrow} \P \left(\frac{\tau_N}{\beta_N} > s\right) \P \left( \frac{\tau_N}{\beta_N} > t \right).
\end{equation*}

Because of the spatial structure, deciding how the typical states shall be defined is a bit more subtle than in the previous complete setting. One important quantity for doing so in the present setting is the density $\rho$ of active neurons under stationarity for the infinite system: $$\rho \triangleq \mu (\{S \subset \Z: 0 \in S\}).$$

Notice that the density is defined with respect to neuron $0$, but could be defined indifferently with respect to any neuron $i \in \Z$ since the dynamics is translation invariant. The typical states for the finite system are the one on which the empirical density of active neurons is not too far below the density of the infinite system. More precisely, we define:
$$ \mathcal{W}_N \triangleq \left\{ W \subset \Z \cap [\m N, N]: \frac{|W \cap [\m N/2,0]|}{N/2+1} > \frac{\rho}{2}, \; \frac{|W \cap [0,N/2]|}{N/2+1} > \frac{\rho}{2} \right\}.$$
Again, by elementary computations one obtains the following upper-bound on the factorization defect:
\begin{multline*}
    \left| \P \left( \frac{\tau_N}{\beta_N} > s + t \right) - \P \left( \frac{\tau_N}{\beta_N} > s \right) \P \left( \frac{\tau_N}{\beta_N} > t \right) \right|\\
    \leq \max_{W \in \mathcal{W}_N} \P \left( \tau_N \neq \tau_N^W \right) 
    + \P \left(\tau_N > \beta_N s, \, \lambda^N(\beta_N s) \notin \mathcal{W}_N\right)
\end{multline*}
The first term in the bound above --- in which $\tau_N^W$ denotes the extinction time of the copy started from $W$, realized together with the full-start process on the common probability space of the Harris construction --- is the one for which one shall identify some renovation pattern on the typical states, while the second term corresponds to the atypical states, which shall be proven unlikely in typical trajectories.

The arguments are a bit more involved than in the previous complete case. The proof --- which makes heavy use of the naturally synchronous coupling allowed by Harris graphical construction, where spikes and leaks are triggered by the same Poisson processes for any copy of the system --- is structurally close to the one established in \cite{cassandro} for the contact process. The coupling allowed by Harris graphical construction is naturally synchronous, since spikes and leaks are triggered by the same Poisson processes for any copy of the system; moreover the infinite version of the system and its non-trivial invariant measure $\mu$ play the role of the proxy to the finite system. The main breaking point with respect to the contact process situation is the absence of self-duality. However, the spiking rates process $(\lambda^N(t))_{t \in \R_+}$ does have a dual, which happens to preserve most of the good properties of the original process.

The configurations in $\mathcal{W}_N$ have two advantages that one can take advantage of to prove that both terms in the bound above vanish: first they are fairly populous configurations, and second they are the ones that the infinite system typically visits in the long run due to the spatial ergodicity of $\mu$. If the infinite system starts from a sufficiently populous configuration --- e.g.\@ $W \in \mathcal{W}_N$ --- then it has probability vanishingly close to $1$ of never reaching the absorbing state, which can be shown by duality. But if this happens, it means that at some point this infinite system will necessarily cross the borders of $\{\m N, \ldots, N\}$. Then, by standard arguments based on the nearest-neighbors nature of the interaction and involving Harris graphical construction, the trajectories of the two finite systems --- starting from $W$ and from the maximal configuration $\{\m N, \ldots, N\}$ respectively --- collide and then stay coupled for eternity. Now, at time $\beta_N s$, either the finite system we're ultimately interested in has reached a configuration in $\mathcal{W}_N$ or it didn't. If it did, then the coupling we just sketched takes care of what happens next since it implies that $\tau_N = \tau_N^W$ with high probability. This addresses the first term in the upper-bound above. The second term is addressed by another coupling between the infinite and finite processes, establishing that if the finite process survives long enough it becomes well described by the restriction of the invariant measure $\mu$ to the state space of the finite process. More precisely, if the system has survived up to time $\beta_N s$, but has not reached $\mathcal{W}_N$, then, due to its tendency to propagate for small $\gamma$, it has probability vanishingly close to $1$ of exceeding spatially the central space $\{\m N/2, \ldots, N/2\}$ --- meaning that there is at least one active neuron both on the left of $\m N/2$ and on the right of $N/2$. But in that case Harris graphical construction provides again a coupling, in which the finite and infinite systems coincide on $\{\m N/2, \ldots, N/2\}$, and monotonicity allows one to conclude that the second term is less than $\mu ({\mathcal{W}_N}^\mathsf{c})$, which vanishes as $N$ diverges by the spatial ergodicity of $\mu$.

Higher dimensions were also considered by \textcite{romaro} --- from a purely computational perspective though --- suggesting that such behavior is not restricted to the one-dimensional $\Z$ with nearest-neighbors interaction, but certainly generalizes to $\Z^d$ with nearest-neighbors interaction. It is quite natural to conjecture that such behavior actually holds in a much more general setting: one shall expect \ref{meta:loss} to hold for any sufficiently connected graph with bounded interaction, even though pinning down the exact conditions under which this happens --- and then establishing a general and rigorous proof --- is certainly not an easy task. 

Thermalization was proven by \textcite{mandre3}, in a formulation slightly different from \ref{meta:therma} though. First, one is usually not so much interested in the statistics on active neurons, but only on the statistics related to the spikes, seen as punctual events, and second, for the sake of simplicity the statement focuses on arithmetic averages. For any finite set $F \subset \Z$ of neurons --- one can think of the neurons which are actually observed in any neuro-physiological experiment, which are usually only a small part of the whole network --- one can define the average number of spikes in $F$ on the interval $[t, t+R]$ for a system of $N$ neurons as: $$\bar{Z}_R^N (t, F) \triangleq \frac{1}{R} \sum_{i \in F} Z^N_i([t,t+R]).$$
When $\gamma < \gamma_c$, one can show that there exists some constant $0 < \rho < 1$, depending only on $\gamma$, representing uniformly the average spiking rate of an individual neuron across time. More precisely \cite{mandre3} shows that, at any time $t>0$: 
$$ \bar{Z}_R^N (t, F) \overset{\P}{\underset{N \rightarrow \infty}{\longrightarrow}} \rho \, |F|.$$

\subsection{Exponential activation function and total/unitary leakage} \label{subsec:metaexponentialactivationfunction}

In \textcite{laxa} the GL model is studied under both total and unitary leakage, in a complete setting --- i.e.\@, $w_{i \rightarrow j}=1$ for all neurons $i\neq j$ --- and with an exponentially increasing activation function, given by $\phi(u)=e^{u}\mathbbm{1}_{u>0}$.
The leakage rate is always equal to one --- in the model with total leakage this means that $\gamma=1$ and in the model with unitary leakage this means that $\alpha(u) = 1$, for any $u \in \R_+$.

We sketch the proof in the total leakage version. The infinitesimal generator for the membrane potentials is hence given by \eqref{eq:GLmodeltotalleakage} and, as usual, $N\geq 2$ denotes the number of neurons, which are indexed by $I_N=\{1,\ldots,N\}$. The time evolution of the system of membrane potentials is denoted by $(U^N(t))_{t \in \R_+}$, taking values in $\mathcal{S}_N = \Z_+^N$.
The state in which all neurons have null membrane potential, denoted $\mathbf{0}_N$, is a trap for the process. Hence, the absorbing region is $\mathcal{A}_N \triangleq \{\mathbf{0}_N\}$ and the metastable region is $\mathcal{M}_N \triangleq \mathcal{S}_N \setminus \mathcal{A}_N$.

The central idea of the proof resides in the identification of an important set of states called \term{ladder lists}, and defined as:
$$
\mathcal{L}_N\triangleq\big\{u \in \mathcal{S}_N : \{u_i : i \in I_N\}=\{0,1,...,N-1\}\big\}.
$$
This set contains the states which have the shape of a ladder from $0$ to $N-1$ with unit step --- once the values of the membrane potentials have been sorted in increasing order --- whence the name. The ladder list set have two basic properties:
\begin{itemize}
\item For any $u,v \in \mathcal{L}_N$, there exists a permutation $\sigma:I_N\to I_N$ such that $u_i=v_{\sigma(i)}$, for any $i\in I_N$.
\item When starting in $\mathcal{L}_N$, if the most likely event occurs --- i.e.\@, the spiking of the neuron with membrane potential equal to $N-1$ --- then the process remains in $\mathcal{L}_N$.
\end{itemize}
Note that the second property above holds due to the complete interaction and the increasing activation function.

As usual the extinction time is defined as $\tau_N \triangleq \inf\{t>0: U^N(t) = \mathbf{0}_N\},$ and \cite{laxa} prove that \ref{meta:loss} holds under $\P_{u}$ for any $u$ in a wide set of sufficiently active initial states --- to be defined more precisely below. The general result for initial configurations in this region of the state space is obtained by proving first \ref{meta:loss} for initial states in the ladder lists, and then extending to the wider regions. 
Thus let  $(l_N)_{N \geq 2}$ be any deterministic sequence of initial states with $l_N \in \mathcal{L}_N$ for all $N \geq 2$.
As usual, the memorylessness property \ref{meta:loss} is established for a macroscopic time defined by $\P_{l_N} ( \tau_N > \beta_N ) = e^{-1}$ which is later replaced by  $ \E_{l_N}(\tau_N)$. Observe that, by the permutation invariance of $\mathcal{L}_N$ and by the symmetry of the interaction, the value of $\beta_N$ is itself invariant with respect to the exact $l_N \in \mathcal{L}_N$ chosen as initial state, which is why it does not appear in the notation. Similarly as in preceding cases, \cite{laxa} proves that, for any pair of positive real numbers $s,t\geq 0$, we indeed have:
\begin{equation} \label{eq:almostladder}
\P_{l_N} \left( \frac{\tau_N}{\beta_N} > s + t \right) \underset{N \rightarrow \infty}{\longrightarrow} \P_{l_N} \left(\frac{\tau_N}{\beta_N} > s\right) \, \P_{l_N} \left( \frac{\tau_N}{\beta_N} > t \right).
\end{equation}

As usual, the proof of the convergence above is based on a distinction between typical and atypical states, the typical states being given by
$$
\mathcal{W}_N \triangleq
\left\{u \in \mathcal{S}_N: \left\{1,\ldots, N - \lfloor \sqrt{N} \rfloor \right\}\subset \{u_a : a \in I_N\}
,\ \bigcap_{a \in I_N} \bigcap_{b \neq a} \{u_a \neq u_b\}
\right\}.
$$
The elements of this set are the lists in which all neurons have different membrane potentials and the membrane potential of the neurons create a ladder of size $N - \lfloor \sqrt{N} \rfloor$. This set generalizes the idea behind the ladder lists $\mathcal{L}_N$. While ladder lists are closed under the occurrence of the most likely event, the set $\mathcal{W}_N$ mimics this property for the occurrence of one of the most likely events: when starting in $\mathcal{W}_N$, if one of the $\lfloor \sqrt{N} \rfloor$ most likely events occurs --- i.e.\@, the spiking of one of the neurons with membrane potential greater than $N - \lfloor \sqrt{N} \rfloor$ --- then the process remains in $\mathcal{W}_N$. Again, this happens thanks to the symmetry induced by the complete  interaction. Moreover, with the activation function that grows exponentially these most likely events actually occur with probability $1$ as $N\to \infty$.
Then, in order to prove \eqref{eq:almostladder} one shows that, for any $N \geq 2$ and $l \in \mathcal{L}_N $, the factorization defect of $\tau_N$ under $\P_l$ obeys the following upper-bound:
\begin{multline*}
\left| \P_l \left( \frac{\tau_N}{\beta_N} > s + t \right) - \P_l \left( \frac{\tau_N}{\beta_N} > s \right) \P_l\left( \frac{\tau_N}{\beta_N} > t \right) \right| \leq \\
\sup_{w \in \mathcal{W}_N}\left\lvert\P_l\left(\frac{\tau_N}{\beta_N}> t\right)-\P_{w}\left(\frac{\tau_N}{\beta_N}> t\right) \right\rvert+ \P_l\left(U^{N}(\beta_Ns)  \notin \mathcal{W}_N,\ \tau_N>\beta_Ns \right).
\end{multline*}

As usual, some kind of renovation argument is used to show that the difference in the first term above vanishes as $N$ diverges. More precisely, consider the synchronous coupling between two copies of the system in which, after each jump of the system (either a spike or a leak), the neurons are ordered --- from smaller to greater membrane potential say --- and paired together with respect to their rank. The synchronicity lies in the fact that the pairing is done in a way that maximizes the probability that neurons at the same position of the ordering jump together. Under such coupling, if at some point the two copies of the system visit the set of ladder lists simultaneously, then the trajectories collide, making the extinction times equal. Furthermore, starting from any two states in $\mathcal{W}_N$ the two copies of the processes visit $\mathcal{L}_N$ together before exiting $\mathcal{W}_N$ --- and therefore before extinction --- with probability tending to one as $N$ diverges. This already proves that the first term of the upper-bound indeed vanishes.

For the second term, similarly as for the complete-interaction model considered in Section \ref{subsec:metabinaryactivationfunction}, one might consider a modification of $(U^{N}(t))_{t \in \R_+}$ where the transition to $\mathbf{0}_N$ has been removed. This modified version again plays the role of a proxy and dominates the original version. Moreover it has a unique invariant measure $\mu = \mu_N$ which, although it cannot be computed explicitly, can be shown to concentrate on the typical states: $\mu(\mathcal{W}_N) \rightarrow 1$ as $N$ diverges. From there one can show that, if there are enough active neurons in the initial configuration --- namely, at least $\sqrt{N}$ of them --- then, conditioned on survival, the probability that the system is in $\mathcal{W}_N$ approaches $1$ as $N$ diverges. This implies that the second term in the upper-bound vanishes as well as $N$ diverges.

Finally, once \eqref{eq:almostladder} has been proven one can extend the result to any initial state in the ``active enough'' set: 
$$\mathcal{V}_N \triangleq \left\{ u \in \mathcal{S}_N: u_k > 0 \text{ for at least } \lfloor \sqrt{N} \rfloor \text{ neurons} \right\}.$$
In order to do so one has to use again the modified immortal process defined above together with the information available on its invariant measure $\mu$ to prove that $\mathcal{L}_N$ is reached almost instantaneously from any state in $\mathcal{V}_N$ as $N$ grows. This together with the ladder list result implies that, for any $u \in \mathcal{V}_N$, we indeed have:
\begin{equation*}
\P_u \left( \frac{\tau_N}{\beta_N} > s + t \right) \underset{N \rightarrow \infty}{\longrightarrow} \P_u \left(\frac{\tau_N}{\beta_N} > s\right) \, \P_u \left( \frac{\tau_N}{\beta_N} > t \right),
\end{equation*}
effectively establishing \ref{meta:loss} for any active enough initial state.

The arguments described above, which deal with the total leakage version of the model, can be adapted to the unitary leakage version with very little additional effort. Furthermore, although thermalization is not discussed in \cite{laxa}, the existence of a unique invariant measure for the process where the transition to $\mathbf{0}_N$ has been removed and the concentration of this invariant measure on $\mathcal{W}_N$ as $N\to \infty$, together with the fact that the process quickly hits $\mathcal{W}_N$ when starting with a sufficiently large number of active neurons, strongly suggest thermalization even though not proven in the precise technical sense of \ref{meta:therma}.

\subsection{Saturated linear activation function and continuous leakage} \label{subsec:metastabilityevamomn}

As mentioned earlier, the continuous leakage version of the system, with finitely many neurons, exponential decay and positive synaptic weights was first studied by \cite{duarteost}. Formally, we consider the system defined by \eqref{eq:GLmodelexp_leakage}, with $|I|=N < \infty$, and $w_{i \rightarrow j} > 0$ for any $i \neq j \in I$. As usual $( U^N(t))_{t \in \R_+}$ denotes the membrane potentials for this finite system. Under rather weak assumptions\footnote{Namely: $\phi(0)=0$, positivity on $\R_+^*$, continuity on $\R_+$ and differentiability at $0$.}  on the activation function, \cite{duarteost} establishes in particular that the system stops spiking in finite time:
$$ \sigma_N \triangleq \sup \left\{ t \geq 0: \Delta Z^N_i(t) > 0 \text{ for some } i \in I \right\} < \infty \text{, almost surely.}$$
The model was then studied by \textcite{evamonm} in a mean-field setting, meaning that:
$$w_{j \rightarrow i}=w/N$$
for any $i \neq j \in I$, for some arbitrary constant $w>0$. In words, whenever a neuron spikes the membrane potential of all the other neurons increases uniformly by an amount which is scaled with respect to the system's volume. As usual a spike resets the spiking neuron's membrane potential to $0$ and, in between spikes, the membrane potential decays in a deterministic fashion --- exponentially fast at a given rate $\alpha>0$. Moreover, we assume that $\phi:\R \rightarrow \R_+$ is a saturated linear function of the following form: $$
\phi(u)= \lambda u \wedge \phi^*.
$$
Here, $\lambda>0$ and $\phi^*>0$ are the parameters of this class of activation functions. The metastable behavior of this system was studied in \cite{evamonm} in the pathwise sense, by establishing in particular that \ref{meta:loss} holds for suitably chosen parameters $\lambda$, $\phi^*$, $w$ and $\alpha$. We omit here the exact definition of the intricate region in the parameter space on which metastability occurs\footnote{Intuitively, large values of $\lambda$ and $w$ keep the system active for a long period of time and large values of $\alpha$ push the system to extinction. Condition $\lambda w>\alpha$ is necessary for the metastability result and the other assumptions on \cite{evamonm} are granted for mainly technical reasons.} and refer to Section 2 of \cite{evamonm}. Instead of proving memorylessness for the last spiking time $\sigma_N$ defined above, they work with a closely related hitting time: namely they consider the first time at which the \term{average spiking rate} is sufficiently small. More precisely, for any $u \in \R_+^N$ define $\bar{\phi}_N(u)$ as $$\bar{\phi}_N(u) \triangleq \frac{1}{N} \sum_{k=1}^N \phi(u_k),$$ and then let $\mathcal{A}_N \triangleq \left\{ u \in \R_+^N : \bar{\phi}_N(u) < \epsilon \right\},$ where $\epsilon$ is a small enough value, namely: $$\epsilon < \frac{\phi^*}{2} \left(1-\frac{\alpha + \phi^*}{\lambda w} \right).$$ Moreover, define the metastable region as the complementary set of $\mathcal{A}_N$, that is $\mathcal{M}_N = \mathcal{A}_N^\mathsf{c}$, and the typical states region $\mathcal{W}_N$ is the same, minus some small gap:
$$\mathcal{W}_N \triangleq \left\{ u \in \R_+^N : \bar{\phi}_N(u) \geq \epsilon + \delta \right\},$$ for some small $\delta > 0$. The extinction time is the hitting time of $\mathcal{A}_N$, $$\tau_N \triangleq \inf \left\{ t \geq 0: U^{N}(t) \in \mathcal{A}_N \right\}.$$
Then \cite{evamonm} proves that there exist positive constants $C$ and $\theta$ such that:
$$ \sup_{t \in \R_+} \sup_{u \in \mathcal{W}_N} \left| \P_u \left(\frac{\tau_N}{\E_u(\tau_N)} > t \right) - e^{-t} \right| \leq C e^{-\theta N},$$
for big enough $N$. Hence, for any $u \in \mathcal{W}_N$, one indeed has \ref{meta:loss} on the scale $\beta_N = \E_u(\tau_N)$, that is: $\tau_N/\beta_N \overset{\mathcal{D}}{\underset{N \rightarrow \infty}{\longrightarrow}} \mathcal{E}(1)$ under $\P_u$. 

The thermalization property \ref{meta:therma} is not established in the literal pathwise sense defined earlier, but the machinery involved in the proof of \ref{meta:loss} is reminiscent of \ref{meta:therma} in a slightly different formalism, belonging to the \emph{McKean--Vlasov mean-field theory}. The general ideas --- which are also based on the study of some limit process and its invariant distribution --- are indeed quite similar. The property of \emph{propagation of chaos} is established. More specifically one proves that, in the $N \rightarrow \infty$ limit, neurons de-correlate and the membrane potential of any of them converges in distribution to a limit process which is described by a non-linear stochastic differential equation given by: \begin{equation}\label{eq:mckean}
  d\bar{U}(t) = -\alpha \bar{U}(t)\,dt + w z_t\,dt
  - \bar{U}(t^-)\int_{\R_+} \mathbbm{1}_{\{z \leq \phi(\bar{U}(t^-))\}}\,\zeta(dt, dz),
\end{equation}
where $\zeta$ is a Poisson random measure on $\R_+^2$ and $z_t \triangleq \E (\phi(\bar{U}(t)))$ represents the mean firing rate at time $t$. The nonlinearity lies in the fact that the dynamics of $(\bar{U}(t))_{t \in \R_+}$ depends on its own law through $(z_t)_{t \in \R_+}$: it is therefore an equation for which the law of the solution is part of the solution. Proposition~4.1 of \cite{evamonm} quantifies the convergence of the particle system to this limit process, showing that the empirical mean $\bar{\phi}_N(U^N(t))$ approximates $z_t$ with error of order $1/\sqrt{N}$.
 
Moreover, in the supercritical regime this mean-field limit has two invariant measures, the trivial one $\delta_0$ that put all mass on the point $0\in \R_+$ and another non-trivial and absolutely continuous invariant measure, say $m$, with positive support on $\R_+$ and which density can be expressed explicitly. The measure $\delta_0$ is unstable while $m$ is globally attracting (for any initial condition other than $\delta_0$) --- see Theorem 2.7, Proposition 2.8 and Theorem 2.12 of \cite{evamonm}. This implies that the Wasserstein distance\footnote{The order-1 Wasserstein distance between two probability measures $p_1$ and $p_2$ on $\mathbb{R}_+$ is defined by $W_1(p_1, p_2) \triangleq \inf_{p} \int |x-y|\,p(dx,dy)$, where the infimum is taken over all probability measures $p$ on $\mathbb{R}_+ \times \mathbb{R}_+$ with marginals $p_1$ and $p_2$.} of two copies of the limit process $(\bar{U}(t))_{t \in \R_+}$, one starting with an initial measure different from $\delta_0$ and the other starting with initial measure $m$, vanishes exponentially fast as time diverges. The finite system with $N$ neurons has of course a trivial invariant measure given by $\nu=\delta_{\mathbf{0}_N} = \delta_0^{\otimes N}$ --- which can be deduced from the preliminary results of \cite{duarteost} alone. Moreover, the convergence of all neurons in the system to the mean-field limit and the fact that the trajectories of this limit process shall evolve as $m$ indeed strongly suggest that, for big $N$, any suitably chosen spatio-temporal statistic on the trajectories of $(U^N(t))_{t \in \R_+}$, taken on a suitably chosen time scale, shall be close to the same statistic as seen from the product measure $\mu = \mu_N \triangleq m^{\otimes N}$, in the sense of \ref{meta:therma}.

To close, let us mention that \cite{evamonm} in fact prove a second metastability statement, which we record here \foreign{en passant} as it pins down the metastable phase more sharply than the result discussed above. Let
$$
p^* \triangleq \int_0^\infty \phi \, dm,
$$
that is, the equilibrium value of the mean firing rate $z_t$ under the non-trivial measure $m$. Whereas the result stated earlier localizes the system away from extinction, one may instead localize it around the equilibrium rate $p^*$ itself. To this end we recycle the notation above, redefining $\mathcal{A}_N$ and $\mathcal{W}_N$ as
$$
\mathcal{A}_N \triangleq \left\{ u \in \R_+^N : \left| \bar{\phi}_N(u) - p^* \right| > \epsilon \right\}
\quad \text{and} \quad
\mathcal{W}_N \triangleq \left\{ u \in \R_+^N : \left| \bar{\phi}_N(u) - p^* \right| \le \epsilon - \delta \right\},
$$
for a small radius $\epsilon > 0$ and a gap $\delta \in (0,\epsilon)$. This time $\mathcal{A}_N$ is the complement of a neighborhood of $p^*$ rather than of the active region, and $\mathcal{W}_N$ the concentric inner band. Considering the hitting time $\tau_N$ of the new $\mathcal{A}_N$, \textcite{evamonm} also establishes that $\tau_N / \E_u(\tau_N)$ likewise converges to a unit-mean exponential as $N$ diverges under $\P_u$, for any $u$ in the new $\mathcal{W}_N$. The sole difference with the previous statement lies in the speed of convergence: in place of the exponential control $Ce^{-\theta N}$, one obtains only a polynomial bound of order $C \ln N / N^{1/4}$. In words, throughout the metastable phase the empirical spiking rate does not merely stay away from extinction: it remains pinned near its mean-field equilibrium value $p^*$, departing from it abruptly at an unpredictable time.

\subsection{A general result for asymptotically exponential exit times} \label{subsec:generalresult}

 The metastability results of \cite{evamonm} described in the previous subsection are actually obtained by applying a general, model-independent theorem proven in the same article --- generalizing the result of \cite{brassesco} on low-noise diffusions --- which we now state on its own. This result is of independent interest: it somewhat isolates in abstract form the proof structure common to all the works surveyed in this section --- as we point out below --- and it continues a line of pathwise, reversibility-free criteria for asymptotically exponential hitting times going back to \cite{metaest1}. Let $(X^N(t))_{t \in \R_+}$ be a time-homogeneous strong Markov process taking values in a Polish space $E$, parametrized by a parameter $N$ representing, in our context, the volume of the system. Consider as usual sets $\emptyset \neq \mathcal{W}_N \subset \mathcal{M}_N \subset E$. For any $A \subset E$ we denote $\tau_{A} \triangleq \inf \{t\geq 0:X^N(t)\in A\}$ and we are specifically interested in the distribution of $\tau_N \triangleq \tau_{\mathcal{A}_N}$ where $\mathcal{A}_N = \mathcal{M}_N^\mathsf{c}$. Consider the three following conditions.
\begin{enumerate}[label=\textbf{C.\arabic*}]
    \item \label{cond:epsilon1} Starting from $\mathcal{W}_N$, the process does not reach $\mathcal{A}_N$ quickly:
    $$
    \sup_{x\in \mathcal{W}_N} \P_x \left(\tau_N \leq s_1\right) \leq \epsilon_1, \text{ for some } \epsilon_1 = \epsilon_1(N) \in [0,1] \text{ and } s_1 = s_1(N) >0.
    $$
    \item \label{cond:epsilon2} Starting from $\mathcal{M}_N \setminus \mathcal{W}_N$, the process quickly reaches either $\mathcal{W}_N$ or $\mathcal{A}_N$:
    $$
    \sup_{x\in \mathcal{M}_N \setminus \mathcal{W}_N} \P_x \left( \tau_{\mathcal{W}_N\cup\mathcal{A}_N} > s_2 \right) \leq \epsilon_2, \text{ for some } \epsilon_2 = \epsilon_2(N) \in [0,1] \text{ and } s_2 = s_2(N) >0.
    $$
    \item \label{cond:epsilon3} The difference between the distribution function of $\tau_N$ for initial states in $\mathcal{W}_N$ is small:
    $$
    \sup_{t \in \R_+}\sup_{x,y\in \mathcal{W}_N} \left| \P_x \left( \tau_N > t \right)-\P_y \left( \tau_N > t \right) \right| \leq \epsilon_3, \text{ for some } \epsilon_3 = \epsilon_3(N) \in [0,1].
    $$
\end{enumerate}

In our context the constants $\epsilon_1$, $\epsilon_2$, $\epsilon_3$, $s_1$ and $s_2$ shall be thought of as functions of the volume $N$ of the system under consideration. Intuitively, the meaning of these constants is as follows. Times $s_1$ and $s_2$ corresponds to time-scales, $s_1$ being meant to be a \emph{long time scale} and $s_2$ a \emph{short time scale}. The epsilons are \emph{failure probabilities} for each condition, meant to be small. While \ref{cond:epsilon1} expresses that typical trajectories dwell in $\mathcal{W}_N$ for a long time (\emph{slow escape}), \ref{cond:epsilon2} that the atypical states $\mathcal{M}_N \setminus \mathcal{W}_N$ are quickly left in favor of either $\mathcal{W}_N$ or $\mathcal{A}_N$ (\emph{fast recurrence}), and \ref{cond:epsilon3} formulates the renovation which takes place within $\mathcal{M}_N$ (\emph{loss of memory of the initial state}).

Assume that the three conditions above hold for epsilons satisfying $\epsilon_1 +\epsilon_2 +\epsilon_3 \leq \tfrac{1}{2}$ and that there exists $x \in \mathcal{W}_N$ such that $\tau_N$ is a.s. finite under $\P_x$. Under these assumptions, one can then prove that there exists $\beta_N >0$ such that $\P_x(\tau_N > \beta_N) \in [1/4,3/4]$. Such a property would of course be trivial for a diffusion, but we're interested in piecewise deterministic Markov processes, for which the law of $\tau_N$ might very well have atoms. The meaning of $\beta_N$ is that of a \emph{typical exit time}, since it lies on a scale more or less of the order of the median (by definition it literally lies inside the inner window of the ``box-plot''). This condition is only there to ensure we're not in the degenerate case where one can't even identify a reasonable candidate for the macroscopic time scale on which metastability is expected to manifest itself.

Granting all assumptions above, \cite{evamonm} proves that $\tau_N/\E(\tau_N)$ is approximately distributed as a unit mean exponential law. Informally, the result proven there states that, additionally to the requisites stated above, if the errors are small --- i.e.\@ $\epsilon_2(N) + \epsilon_3(N) \approx 0$ --- and the short time scale (the recurrence time) is negligible with respect to the typical exit time --- i.e.\@ $s_2(N)/\beta_N \approx 0$ --- then, for any initial distribution $\nu$ on $\mathcal{M}_N$ that isn't about to exit immediately --- i.e.\@ $\P_\nu (\tau_N \leq s_2(N)) \approx 0$ --- one has $  \P_\nu(\tau_N/\E_\nu (\tau_N) > t) \approx e^{-t}$ uniformly over $t$. The approximate exponentiality of the hitting time is stated together with an explicit control over the quality of this approximation, and the usefulness of such a result with respect to our metastable concerns of course depends crucially on the control provided.

Formally, fix some $C,\delta > 0$ and then choose the biggest $M = M(N) >0$ such that
$$
(s_2(N)/\beta_N) \vee(\epsilon_2(N)+\epsilon_3(N)) \leq Ce^{-\delta M(N)}.
$$
Above $C$ and $\delta$ are absolute constants while $M$ shall be thought of as an increasing function of $N$ in our context. The first statement of Theorem 5.3 in \cite{evamonm} states that there exists $M_0 = M_0(C,\delta)$ such that, if $M$ can be chosen so that $M(N) \geq M_0$, then $\sup_{y \in \mathcal{M}_N}\E_y(\tau_N)<+\infty$. The second statement moreover asserts that, if $\nu$ is a measure on $\mathcal{M}_N$ satisfying $\P_\nu(\tau_N\leq s_2(N))\leq Ke^{-\vartheta M(N)}$ for some arbitrary positive constants $K$ and $\vartheta$, then there exist two constants $C_1=C_1(\delta,C,K,\vartheta)>0$ and $C_2=C_2(\delta,\vartheta)>0$, which depend neither on $N$ nor on the choice of $M$, such that:
\begin{equation} \label{eq:generalmeta}
    \sup_{t \in \R_+} \left|\P_\nu \left(\frac{\tau_N}{\E_\nu(\tau_N)}>t \right)-e^{-t} \right| \leq C_1M(N)^3e^{-C_2 M(N)}.
\end{equation}

Observe that, so far, $\epsilon_1$ and $s_1$ play quite a minor role, since they only appear in the condition $\epsilon_1 + \epsilon_2 + \epsilon_3 \leq \tfrac{1}{2}$, which only serves as a guardrail to discard the degenerate case and should anyway be easy to satisfy since $\epsilon_2 + \epsilon_3$ is later required to vanish exponentially fast. In particular $\epsilon_1$ and $s_1$ do not seem to play any role in the convergence rate of \eqref{eq:generalmeta}. Nevertheless, this is only so if one is able to actually establish the two other conditions independently of $\epsilon_1$ and $s_1$, and in particular to prove \ref{cond:epsilon3} directly, which is usually difficult to do. In practice, \eqref{eq:generalmeta} is established via a sufficient condition provided by \cite{evamonm} and formulated in terms of a coupling. Namely, grant \ref{cond:epsilon1} for some $\epsilon_1$ and $s_1$, and assume moreover that there exists an $\epsilon_4 = \epsilon_4(N) \in [0,1]$ and a copy $(Y^N_t)_{t \in \R_+}$ of $(X^N_t)_{t \in \R_+}$ such that for all $x,y \in \mathcal{W}_N$ one has $\P_{x,y} (X^N(t) = Y^N(t) \text{ for all } t \geq s_1(N)) \geq 1-\epsilon_4(N)$, where $\P_{x,y}$ denotes the joint law under which $(X^N(0),Y^N(0))=(x,y)$. Then \ref{cond:epsilon3} holds with $\epsilon_3 = 2 \epsilon_1 + \epsilon_4$. Hence, for all practical purposes, the condition on $\epsilon_1$ will actually be as stringent as for the other two, for $\epsilon_1$ is actually embedded in $\epsilon_3$.

All this means is that, if one is interested in proving \ref{meta:loss}, then it is enough to prove \ref{cond:epsilon1}--\ref{cond:epsilon3} with sufficiently good $\epsilon_1$, $\epsilon_2$, $\epsilon_3$, $s_1$ and $s_2$. Here ``sufficiently good'' means that, as $N$ diverges, the epsilons shall vanish fast enough while $s_2$ shall not grow too fast --- where the exact meaning of ``fast enough'' and ``not too fast'' depends on how fast the targeted $M$ grows. Concerning the saturated-linear model of the previous section, the metastability result of \cite{evamonm} is recovered by a direct application of this general result on the potential process $(U^N(t))_{t \in \R_+}$. Conditions \ref{cond:epsilon1} and \ref{cond:epsilon2} follow from propagation of chaos and the properties of the limit process described in previous section, together with the study of an auxiliary Markov process that lower-bounds the evolution of the total jump rate and for which a large deviation principle holds; this is what makes the relevant hitting times exponentially large in $N$. Condition \ref{cond:epsilon3} follows from Proposition 5.8 of \cite{evamonm}, establishing the relevant coupling for the sufficient condition to hold.

It is worth explaining why the two domains of the previous subsection yield qualitatively different rates, $Ce^{-\theta N}$ in the first case and $C\ln N/N^{1/4}$ in the second. In both, the rate is governed by how large the scale-separation parameter $M$ may be taken, which is in turn limited by how fast the failure probabilities $\epsilon_2 + \epsilon_3$ vanish --- the smaller they are, the larger $M$, and the sharper the bound \eqref{eq:generalmeta}. In the first case $\mathcal{W}_N = \{\bar\phi_N \ge \epsilon+\delta\}$ is a region of high global activity, and leaving it is a \emph{large-deviation} event for the auxiliary jump-rate process; all the failure probabilities are then exponentially small in $N$, so one may take $M(N)=N$ and the bound collapses to an exponentially small error. In the second case $\mathcal{W}_N = \{|\bar\phi_N-p^*|\le\epsilon-\delta\}$ is a narrow window around the equilibrium value $p^*$, within which the loss of memory \ref{cond:epsilon3} is no longer a large-deviation phenomenon but a matter of \emph{propagation of chaos} --- and such fluctuation estimates are only polynomial. Through the reduction $\epsilon_3 = 2\epsilon_1 + \epsilon_4$, the slow-escape error $\epsilon_1$, of order $N^{-1/4}$, dominates $\epsilon_3$; hence $\epsilon_2 + \epsilon_3$ vanishes only polynomially and $M$ can be taken no larger than of order $\ln N$, which yields only a polynomial rate. This illustrates, a posteriori, the role of $\epsilon_1$: exponentially small and harmless in the first case, it is precisely the bottleneck in the second, entering $\epsilon_3$ through the coupling and fixing the final rate.

\subsection{Related works} \label{sec:related}

Below we sum up some metastability results concerned with systems which are either inspired by the GL model but do not completely fit the framework introduced in Section \ref{sec:GLModels} or studied via numerical simulations rather than rigorous proofs of metastability in the pathwise standard.

\subsubsection*{A discrete-time version of the GL model}

The model studied in \textcite{ludmiguel} is a discrete-time version of the GL model in which, unlike the original one of \cite{glmodel}, at most one neuron spikes at a time. At each step a uniformly chosen neuron with potential $x \in \R_+$ spikes with probability $\phi(x) = x \wedge 1$: it then resets to $0$ while every other potential is updated by $y \mapsto \alpha y+w$. If no neuron spikes, all potentials simply leak, that is: $y \mapsto \alpha y$. Here $\alpha \in (0,1)$ is the leakage rate and $w>0$ the synaptic weight. Formally, for a system with $N\geq 2$ neurons, the spiking maps $\pi_i$ are defined for any $1 \leq i \leq N$ and $u \in \R_+^N$ by
\begin{equation*} 
\left(\pi_i(u) \right)_j \triangleq
\begin{cases}
\alpha u_j+w, &\text{ if } j\neq i, \\
0 &\text{ if } j= i,
\end{cases}
\end{equation*}
and the leakage map $\pi^{\dagger}$ is defined by
\begin{equation*} 
\left(\pi^{\dagger}(u) \right)_j \triangleq
\alpha u_j.
\end{equation*}
The evolution of the membrane potentials $(U(t))_{t \in \N}=(U_i(t):i\in I)_{t \in \N}$ is a Markov chain with transition probabilities:
$$
\P \left( U(n)=\pi_i(u) \mid U(n-1)=u \right)=\frac{\phi(u_i) }{N},
$$
and
$$
\P \left( U(n)=\pi^{\dagger}(u) \mid U(n-1)=u \right)=1-\sum_{i \in I}\frac{\phi(u_i)}{N}.
$$

This system stops spiking a.s.\@ in finite time for any $w>0$ and $\alpha \in (0,1)$. The mean extinction time grows with the network size only when $u^\ast > 1$, where $u^\ast \triangleq w/(1-\alpha)$. In this regime the mean potential exhibits metastable behavior, fluctuating around the positive value $u^\ast$ for most of the time before extinction; more precisely, \cite{ludmiguel} shows that the ordered potentials then remain close to the set
$$
\mathcal{W}_N=\{u \in \R_+^N: 0=u_1 \leq u_2 \leq \ldots \leq u_N, (1-\alpha^{i-1})u^\ast \leq u_i \leq u^\ast\}.
$$
The authors also provide bounds on the expected last spiking time, and illustrate the memorylessness property \ref{meta:loss} numerically, without rigorously proving metastability in the strong pathwise sense.

\subsubsection*{Replica-mean-field limits}

\textcite{taille} studies metastability for the GL model through the replica-mean-field (RMF) limit. Unlike the mean-field limits considered in Section~\ref{subsec:metastabilityevamomn}, where synaptic weights are scaled as $w/N$ and the finite-size fluctuations --- which drive the metastable transitions --- are washed out as $N\to\infty$, the RMF limit is a mean-field limit that \emph{preserves} finite-size effects.

Consider the continuous-time GL model with exponential leakage, and assume the activation function is allowed to vary from neuron to neuron, with the general exponential form
$$
\phi_i(u)=h_i\,e^{\kappa_i u},
$$
with $h_i,\kappa_i>0,\ i\in I$. Two features distinguish this choice from the activation functions considered in the previous sections. First, $\phi_i(0)=h_i>0$ and hence each neuron retains a nonzero spontaneous rate, so that there is here \emph{no quiescent absorbing state} --- in contrast with every setting above, where $\phi(0)=0$ makes $\mathbf{0}_N$ an absorbing trap. Second, because $\phi_i$ is positive on the whole real line, negative weights $w_{i \rightarrow j}<0$ are allowed, the membrane potentials range over $\R$ rather than $\R_+$ without ever violating the non-negativity of the intensity. This makes \cite{taille} the only work reviewed here accommodating for inhibition.

Given such a network of $N\geq 2$ neurons, the RMF approach constructs an enlarged network made of $R\geq 2$ replicas of the original one, with neuron set $\{(i,r):i\in I,\ r\in\{1,\ldots,R\}\}$. At each spiking time of a neuron $(i,r)$, its membrane potential is reset to zero and, for each $j\in I\setminus\{i\}$, a replica $r'\in\{1,\ldots,R\}$ is drawn uniformly at random and $w_{i \rightarrow j}$ is added to the membrane potential of $(j,r')$; these draws are independent across spiking times and across neurons. Between spikes, the membrane potential of each $(i,r)$ relaxes to $0$ at exponential speed with rate $\alpha_i$.

The point of this randomized routing is that, as $R\to+\infty$, the interaction between any two fixed neurons occurs with vanishing probability, so that within a representative replica the neurons become asymptotically independent and each is driven by independent Poisson processes --- the so-called Poisson hypothesis. The resulting limit dynamics is therefore parametrized by stationary rates that can be computed by numerical methods and leveraged to study metastability.

The network analyzed in \cite{taille} consists of $N=40$ neurons arranged in two symmetric groups. Each group comprises a cluster of excitatory neurons and a dedicated cluster of inhibitory neurons ($10$ neurons per cluster). Excitatory neurons excite their own cluster, while the cross-inhibition between the two groups is mediated by the inhibitory clusters, endowing the circuit with strong mutual inhibition. A computational study shows that when the excitatory synaptic weight is small, the high level of inhibition keeps the whole network at low activity; when it is large, the network switches between periods of high activity of group $1$ together with low activity of group $2$, and periods of the reverse configuration. These alternations in the global activity are apparently metastable (see Figure~10 of \cite{taille}).

Studying the RMF limit of this system makes it possible to distinguish the two regimes: for small excitatory weight the limit system admits a single stable solution, whereas in the bistable regime it admits two. Beyond detecting multistability, the limit system yields estimates of the mean membrane potential, the mean activity of each neuron, and the higher moments of the membrane potential in each phase.

\subsubsection*{A social network model}

By interpreting membrane potential as social pressure, \textcite{galveslaxa} introduces an opinion model inspired by the GL framework and proves metastability using similar ideas to those in \cite{laxa}. Consider a finite set of social actors $I_N=\{1,2,...,N\}$, with $N \geq 3$. Associated to each actor, we have a point process. We consider marks ($\m 1$ or $\p 1$) on the point processes, indicating the successive moments in which a social actor expresses a ``favorable'' ($\p 1$) or ``contrary'' ($\m 1$) opinion.  Therefore, this model is a system of interacting marked point processes with memory of variable length that describes the temporal evolution of a social network.

The orientation and rate at which an actor expresses an opinion are influenced by the social pressure exerted on it. The social pressure of an actor is reset to $0$ when the actor expresses an opinion, and, simultaneously, the social pressures on all other actors change by one unit in the direction of the expressed opinion. The rate at which actors express opinions and the tendency of each actor to express an opinion in the direction of its social pressure increase exponentially with the \term{polarization coefficient} $\gamma\geq 0$.
This model is formally described as follows.

A state of the system is a vector of \term{social pressures} which are roughly analogous to the membrane potentials of the neurons in GL models. Social pressures are integer-valued, and the state space is thus $\mathcal{S}_N \triangleq \Z^N$. To describe the time evolution of the social network, we introduce a family of maps on $\mathcal{S}_N$. For any actor $a \in I_N$, for any opinion $o \in \{\m 1,\p 1\}$ and for any current vector of social pressures $u \in \mathcal{S}_N$, 
we define the new vector $\pi_{a,o} (u)$ as follows. For all $b \in I_N$, 
$$
\left(\pi_{a,o}(u)\right)_b\triangleq
\begin{cases}
u_b+o &\text{ if } b\neq a, \\
0 &\text{ if } b  = a.
\end{cases}
$$
For a polarization coefficient $\gamma\geq 0$, the time evolution of the social pressures $(U^{\gamma} (t))_{t \in \R_+}$ is a Markov jump process taking values in $\mathcal{S}_N$ and with infinitesimal generator given by
\begin{equation*} 
\mathcal{G}f(u)= \sum_{o \in \{\m 1,\p 1\}}\sum_{b\in I_N}
e^{o \gamma u_b}\left(f(\pi_{b,o}(u))-f(u)\right),
\end{equation*}
for any  bounded function $f:\mathcal{S}_N  \to \mathbb{R}$. Here the relevant limit parameter is the polarization coefficient $\gamma$ rather than the volume $N$ (which is fixed), so it is $\gamma$ that we carry as a superscript on the process; as before, the law and expectation of the process started from $u \in \mathcal{S}_N$ are written $\P_u$ and $\E_u$. 

In this system, the changes between consensus on the network display a metastable behavior. By \term{consensus}, we mean configurations in which the social pressures of all actors push in the same direction, that is, any configuration in one of the two following sets:
$$
\mathcal{C}_N^+ \triangleq \left\{u \in \mathcal{S}_N \setminus \{\mathbf{0}\}: u_a \geq 0 \ \forall a \in I_N \right\} \quad \text{and} \quad \mathcal{C}_N^- \triangleq \left\{u \in \mathcal{S}_N \setminus \{\mathbf{0}\}:u_a \leq 0 \ \forall a \in I_N \right\}.
$$
Formally \cite{galveslaxa} proves that \ref{meta:loss} holds in the $\gamma \to \infty$ limit, when one considers the hitting time of the negative consensus $\mathcal{C}_N^-$ when starting from an arbitrary state in the positive consensus $\mathcal{C}_N^+$. By symmetry this is of course equivalent to \ref{meta:loss} holding for the hitting time of $\mathcal{C}_N^+$ with an initial state in $\mathcal{C}_N^-$. 

In order to prove such a result, one first considers particular subsets $\mathcal{L}_N^+$ and $\mathcal{L}_N^-$ of $\mathcal{C}_N^+$ and $\mathcal{C}_N^-$ respectively, called ladder lists, and defined as:
\begin{align*}
    &\mathcal{L}_N^+\triangleq \left\{u \in \mathcal{S}_N: \{u_a : a \in I_N\}= \left\{0, 1, \ldots, N - 1 \right\} \right\}\\
    \text{and } &\mathcal{L}_N^- \triangleq \left\{u \in \mathcal{S}_N: \{u_a : a \in I_N\}=\left\{0,\m 1,\ldots, \m(N - 1)\right\}\right\}.
\end{align*}
For any configuration $l \in \mathcal{L}_N^+$, property \ref{meta:loss} indeed holds as $\gamma \to \infty$, under $\P_l$, for the following hitting time:
$$
\tau_{\gamma}\triangleq \inf\{t>0: U^{\gamma}(t) \in \mathcal{L}_N^-\}.
$$
As usual, the memorylessness property \ref{meta:loss} is established for a macroscopic time $\beta_\gamma$, defined by $\P_l ( \tau_{\gamma} > \beta_\gamma ) = e^{-1}$ --- which is later replaced by  $ \E_l (\tau_{\gamma})$ --- by showing that, for any pair of positive real numbers $s,t\geq 0$, the following holds
\begin{equation*}
 \P_l\left(\frac{\tau_{\gamma}}{\beta_{\gamma}}>s+t \right)\underset{\gamma \rightarrow \infty}{\longrightarrow}\P_l\left(\frac{\tau_{\gamma}}{\beta_{\gamma}}>s \right) \, \P_l\left(\frac{\tau_{\gamma}}{\beta_{\gamma}}>t \right).    
\end{equation*}
Similarly as in the results reviewed in Section \ref{sec:GLModels} the proof relies on a distinction between typical and atypical states, the set of ladder lists $\mathcal{W}_N\triangleq\mathcal{L}_N^+\cup \mathcal{L}_N^-$ being the typical states in that case. Just as in Section~\ref{subsec:metaexponentialactivationfunction}, the relevance of ladder lists stems from their closure under the most likely transition: starting from a positive ladder, the highest-pressure actor is by far the most prone to speak (its rate $e^{\gamma(N-1)}$ dominating all others as $\gamma \to \infty$), and when it does, it resets to $0$ while every other pressure shifts up by one, so that the system does indeed stay in $\mathcal{L}_N^+$. The symmetry induced by the complete-interaction implies that the law of $\tau_{\gamma}$ under $\P_u$ is the same for any $u \in \mathcal{L}_N^+$. One can take advantage of such symmetry in order to establish that, for any $\gamma\geq 0$ and $l \in \mathcal{L}_N^+$, the following upper-bound holds: 
\begin{equation*}
\left\lvert\P_l\left(\frac{\tau_{\gamma}}{\beta_{\gamma}}>s+t \right)-\P_l\left(\frac{\tau_{\gamma}}{\beta_{\gamma}}>s \right)\P_l\left(\frac{\tau_{\gamma}}{\beta_{\gamma}}>t \right) \right\rvert\leq 
\P_l\left(U^{\gamma} (\beta_{\gamma}s) \in \mathcal{S}_N\setminus\mathcal{W}_N\right).
\end{equation*}
Then one proves that the process $(U^{\gamma}(t))_{t \in \R_+}$ has a unique invariant measure $\mu_\gamma$, and that $\mu_\gamma$ concentrates on $\mathcal{W}_N$ as $\gamma$ diverges, while the hitting time of $\mathcal{W}_N$ becomes vanishingly small for any non-null initial configuration. These two ingredients allow simultaneously to conclude that the right-hand term above vanishes as $\gamma \to \infty$, and then to generalize \ref{meta:loss} to the consensus switching case discussed initially.
Property \ref{meta:therma} is not proven for this model, even though the instantaneous hitting time on the ladder lists on which $\mu_\gamma$ concentrates strongly suggests that such thermalization should hold.

A mean-field version of the model presented above has been studied in \textcite{evalaxa} under much more general conditions on the jump rates. This last work does not study metastability though, but proves instead propagation of chaos: when $N$ diverges each actor's individual social pressure can be described by a limit equation which plays in the present setting a role analogous to \eqref{eq:mckean} in Section \ref{subsec:metastabilityevamomn}. For a suitable choice of parameters \cite{evalaxa} manages to prove that the solution of this mean-field limit equation possesses two invariant probability measures, one supported on the positive real numbers and the other supported on the negative real numbers. The analogies with the situation of \cite{evamonm} suggest that one should be able to obtain similar metastability results by adapting the ideas developed there.

\subsubsection*{A model with synaptic facilitation}

\textcite{pouzatandre} introduces a variant of the GL model incorporating short-term synaptic facilitation, motivated by experimental observations of sustained activity in prefrontal cortex networks during delayed-response tasks. The authors argue that the kind of transient information storage required by working memory naturally fits into the conceptual framework of metastability, and propose a minimal stochastic model designed to exhibit such a phenomenon while remaining amenable to analytical and numerical
treatment.

The model can be seen as a GL system without leakage, with hard threshold activation function and complete graph of interactions, augmented with an additional binary variable per neuron representing the facilitation state of its synapse. Specifically, to each neuron $i \in \{1,\ldots,N\}$ is associated, besides the usual membrane potential process $(U^N_i(t))_{t \in \R_+}$, a facilitation process $(F^N_i(t))_{t \in \R_+}$ taking values in $\{0, 1\}$. The synapse of neuron $i$ is said to be \term{facilitated} at time $t$ whenever $F^N_i(t) = 1$. The activation function is a binary function reminiscent of \cite{ferrari}, in the sense that neurons may spike only when their membrane potential has reached a fixed threshold $\theta \in \Z_+^*$, in which case they do so at constant rate $\kappa > 0$. In other words, the activation function reads $$\phi(u) = \kappa \, \mathbbm{1}_{u \geq \theta},$$ and the membrane potential of any neuron can be assumed without loss of generality to take value in the finite set $\{0, 1, \dots, \theta\}$. The synaptic weights are uniform, $w_{i \rightarrow j} = 1$ for all $i \neq j$.

The defining feature of the model is the role played by the facilitation state in the transmission of spikes, leading to a refinement of the standard spiking rule. When a neuron $i$ emits a spike at some time $t$, its membrane potential is reset to $0$ and its synapse becomes (or remains) facilitated --- that is $F^N_i(t^+) = 1$. However, the spike is transmitted to the rest of the network only if the synapse was \emph{already} facilitated at the moment of the spike: if $F^N_i(t^-) = 1$ the spike is called \term{efficient}, and the membrane potential of every other neuron $j \neq i$ which is below threshold increases by one unit; if $F^N_i(t^-) = 0$ on the contrary, the spike is called \term{inefficient} and has no effect on the rest of the network. Furthermore, the facilitation of each neuron is lost spontaneously at rate $\lambda > 0$, regardless of its membrane potential. The interplay between this spike-induced facilitation and its spontaneous decay is the mechanism by which the model captures the phenomenon of short-term synaptic facilitation.

In a similar fashion to the binary system with complete interaction discussed in Section~\ref{subsec:metabinaryactivationfunction} --- and because the dynamics is invariant under permutations of the neurons --- the analysis is conveniently carried out at the level of the aggregated process $(S^N(t))_{t \in \R_+}$ counting, for each pair $(i,j) \in \{0, \dots, \theta\} \times \{0,1\}$, the number $S^N_{i,j}(t)$ of neurons with
membrane potential $i$ and facilitation state $j$. The process $(S^N(t))_{t \in \R_+}$ is a continuous-time Markov chain on the finite state space $$\mathcal{S}_N \triangleq \left\{s \in \{0,\ldots,N\}^{\{0, \ldots,\theta\} \times \{0,1\}}: \sum_{i,j} s_{i,j} = N \right\}.$$
An immediate consequence of the rules described above is that the configurations in which all neurons are unfacilitated and have sub-threshold membrane potential form a trapping region for the system, since neither a spike nor a facilitation loss can occur from there. The authors give a wider characterization of the absorbing region $\mathcal{A}_N \subset \mathcal{S}_N$ as the set of configurations from which the trajectory is forced to reach the trapping region in a bounded number of transitions --- because there is not enough activity or because the activity is too concentrated on neurons with unfacilitated synapse. Such states can be characterized explicitly as follows
$$\mathcal{A}_N \triangleq {\bigcup_{k=1}^\theta \left\{ s \in \mathcal{S}_N: \sum_{i=k}^\theta s_{i,1} \leq \theta - k \right\}} \cup \left\{ s \in \mathcal{S}_N: s_{\theta,0} + \sum_{i=0}^\theta s_{i,1} \leq \theta \right\}.$$
In order to define the metastable region one needs to take out another class $\mathcal{T}_N$ of transient states defined by $\mathcal{T}_N \triangleq \left\{ s \in \mathcal{S}_N: s_{i,0} + s_{i,1} = 0 \text{ for some } 0 \leq i \leq \theta   \right\}$. Then one lets $$\mathcal{M}_N \triangleq \mathcal{S}_N \setminus \left(\mathcal{A}_N \cup \mathcal{T}_N \right).$$

A distinctive contribution of \cite{pouzatandre} is the use of the theory of \emph{quasi-stationary distributions} (QSD) to investigate the metastable behavior of the system. In order to do so one first shows that $\mathcal{M}_N$ is an irreducibility class for $(S^N(t))_{t \in \R_+}$. This property, combined with the finiteness of the state space, implies by classical results the existence and uniqueness of a probability measure $\mu$ supported by $\mathcal{M}_N$, called the quasi-stationary distribution of $(S^N(t))_{t \in \R_+}$, and characterized by $$\P_\mu\big( S^{N}(t) = z \mid \tau_N > t \big) = \mu(z).$$
Here $\tau_N$ denotes as usual the absorption time of the system, that is the hitting time of $\mathcal{A}_N$. Another well-known fact is that, when the system starts from the QSD, the absorption time $\tau_N$ is \emph{exactly} exponentially distributed: $\P_\mu(\tau_N > t) = e^{-\eta_N t}$, with rate $\eta_N > 0$ given by the opposite of the Perron--Frobenius eigenvalue of the sub-generator obtained by restricting the infinitesimal generator of $(S^N(t))_{t \in \R_+}$ to $\mathcal{M}_N$.

This QSD perspective bridges the gap with the pathwise approach as follows. If, as $N$ diverges, the mixing time toward $\mu$ --- that is, the time required for the conditioned process to come sufficiently close to its QSD --- is asymptotically negligible with respect to the absorption time, then $\tau_N / \E_x(\tau_N)$ obviously converges in distribution to a unit-mean
exponential random variable for any initial state $x \in \mathcal{M}_N$, recovering \ref{meta:loss} on the macroscopic time scale $\beta_N = \eta_N^{-1}$. Establishing this rigorously is left as an open problem, but the authors provide substantial numerical evidence: a direct computation of $\mu$ via the Perron--Frobenius spectral analysis discussed above is carried out for small
$N$, and indeed simulations show that the empirical distribution of the surviving trajectories converges to $\mu$ on a time scale much shorter than $\beta_N$. Thermalization in the sense of \ref{meta:therma} is strongly suggested by these numerical experiments. Moreover the QSD approach gives an alternative formulation of thermalization in which the approximate macroscopic description via the invariant distribution of the infinite system is replaced by an exact description (modulo the relaxation delay) which is valid at the mesoscopic scale (i.e.\@ $N$ need not be very big). Even though there is usually no closed-form formula for the QSD, it can give interesting numerical insight by simple linear algebra and spectral analysis.

\section{Open Problems} \label{sec:conclusion}

We close our review with a list of open questions and further developments regarding metastability in the GL framework which have yet to be addressed.

\vspace{0.2 cm}

\noindent\textbf{Beyond the mean-field setting.} All the metastability results presented above were obtained in a mean-field or complete-interaction setting --- that is, in a perfectly symmetrical world in which the interaction has no spatial or topological structure whatsoever --- with the only exception of the one-dimensional lattice results obtained in \cite{mandre1,mandre3}, which were nonetheless obtained at the expense of drastic simplification with respect to the remaining degrees of freedom of the model (binary neurons with total leakage). Generalization of \cite{mandre1} to the $d$-dimensional lattice --- which is known to be a biologically relevant network structure in some very specific areas of the brain, such as the cerebellar cortex --- have been investigated numerically in \cite{romaro}. Even if no rigorous proof is available for the general $d$-dimensional lattice case, such proof is almost certainly not out of reach of modern technology --- adapting ideas already available in the interacting particle systems literature. However, extension of these results to more realistic versions of the GL model and more interesting, spatially structured, interaction graphs (trees, power-law graphs, Erdős–Rényi...), is at the same time very relevant and quite challenging, and requires the development of new ideas.

We shall mention that lattice results obtained in \cite{mandre1,mandre3} crucially depend on the phase transition established in \cite{ferrari} for the infinite network model. A similar phase transition has already been obtained in networks defined on homogeneous trees. In particular, \cite{amarcos} studies the GL model on homogeneous trees of degree $d\geq 2$, with binary neurons and unitary leakage. In this setting it is proven that the network can become extinct, survive locally, or survive globally depending on the value of the leakage rate. The existence of a regime in which the system globally survives suggests that \ref{meta:loss} might hold for a sequence of finite networks that increasingly exhausts the homogeneous tree of a given degree.

\vspace{0.2 cm}

\noindent\textbf{What about inhibition?} As we've seen, the GL model is well-defined under very weak constraints on the weights $w_{i \rightarrow j}$, requiring only $\sup_{i \in I}\sum_{j \in I}|w_{j \rightarrow i}| < +\infty$, so nothing prevents introducing inhibitory neurons via negative weights. Yet all the rigorous pathwise results above assume a purely excitatory setting --- at odds with the estimate that inhibitory neurons make up roughly 20\% of a typical biological network. The only work considering inhibition is \cite{taille}, and it relies mostly on numerical methods. The technical obstruction in the lattice results of \cite{mandre1, mandre3} resides in the loss of attractiveness --- which underpins the Harris coupling, the finite/infinite comparison, and the spatial-ergodicity argument. In the complete-interaction settings of \cite{mandre2} and \cite{laxa} the introduction of inhibitory neurons breaks the permutation symmetry, which allowed reduction to an aggregated process in one case, and was behind the relevance of ladder lists in the other case. Finally, in \cite{evamonm} the result of metastability depends on a synchronous coupling in which two copies of the system contract because the excitatory drift forces both into the saturated regime where their rates coincide. Negative weights break decisive elements of the proof in each of the cases above. Moreover \textcite{taille} formulates metastability in terms of a switch between two active phases rather than escaping toward extinction --- suggesting inhibition is bound up with the multistability discussed next, making it a central theme to be investigated.

\vspace{0.2 cm}

\noindent\textbf{Multistability.} In the perspective adopted in the works discussed in preceding sections, one usually proves the existence of a single metastable equilibrium in which the system remains during a long and unpredictable time before switching to the quiescent equilibrium. In biological terms, this view only makes sense if one idealizes the system considered as a small sub-network of a bigger system, which has been artificially isolated from it. When they use the expression ``metastable brain'', neuroscientists usually refer to phenomena of multistability, where multiple metastable equilibria are visited in sequence.  Big neural networks (e.g.\@ the entire brain) are usually composed of multiple (functionally consistent) sub-networks, which is the reason why metastability is such an important feature: it allows a balance between adaptability and stability, by regularly switching between different functional sub-networks while avoiding a completely chaotic out-of-equilibrium behavior. It is actually quite natural to conjecture that two important mechanisms behind the occurrence of such multistable dynamics are the complexity of the interaction graph and the delicate intertwining between excitation and inhibition, which would explain why these are not observed in the settings discussed in this review.

In view of that, an appealing definition of multistability could be the following. Denote by $\mathcal{M}_1, \ldots, \mathcal{M}_d \subset \mathcal{S}$ the $d$ metastable regions, and for any $S \subset \mathcal{S}$ let $\tau_S$ be the hitting time of the region $S$. Then the system is multistable if the following holds:
$$ \frac{\sup_{x \not\in \cup_j \mathcal{M}_j } \E_x\left( \tau_{\cup_j \mathcal{M}_j} \right)}{\inf_i \inf_{y \in \mathcal{M}_i} \E _y\left( \tau_{\cup_{j \neq i} \mathcal{M}_j} \right)} \longrightarrow 0,$$
where the limit is taken as the number $N$ of neurons in the system diverges. This definition encapsulates nicely the multiplicity of time scales typical of metastability, and, for a Markovian system, it is expected to entail that the transitions from one region to the other are exponentially distributed. A recent and quite powerful framework for working out such property is the potential theoretic approach developed by A. Bovier and co-authors --- see \cite{bovierholl} for an in-depth account on the subject. Nonetheless much of this theory has been developed only for reversible processes, a property that the systems we consider do not satisfy. Formulations better suited to non-reversible dynamics have nonetheless been proposed, notably the characterization of tunneling for continuous-time Markov chains and the martingale approach of Beltrán and Landim \cite{beltranlandimcharacterization, martingaleapproach}, as well as the more recent resolvent approach \cite{resolventapproach}; we refer to the survey \cite{metastabilitybooklandim} for an overview. Adapting these tools to the multistable, non-reversible regime of GL models is an appealing direction. 

\vspace{0.2 cm}

\noindent \textbf{Thermalization.}
Besides the standard time-scale, the pathwise approach characterizes metastability by the existence of two additional time-scales: a macroscopic one, on which one might observe the tunneling outside the metastable region described by \ref{meta:loss} and the mesoscopic one, on which one observes the thermalization described by \ref{meta:therma}. Most of the works that study metastability using the pathwise approach focus on \ref{meta:loss}, without taking into account \ref{meta:therma}. For the GL model, the only exception is \cite{mandre3}, in which thermalization is proven in a formulation slightly different from \ref{meta:therma}, considering the model in a one-dimensional lattice with binary neurons and total leakage (see Section \ref{subsec:metabinaryactivationfunction}). 
Although the other metastability results presented in Section \ref{sec:metastability_GLmodel} come together with results that suggest that a thermalization property similar to \ref{meta:therma} holds, the work done in \cite{mandre3} illustrates that even after proving \ref{meta:loss}, to prove \ref{meta:therma} it is necessary to deal with non-trivial technical difficulties.

\vspace{0.2 cm}

\noindent \textbf{A general metastability result.}
\textcite{evamonm} prove a general result pinning down three conditions under which a general Markov process has an approximately exponentially distributed exit time --- which we reviewed in Section \ref{subsec:generalresult} above. The authors then use this general result to prove that \ref{meta:loss} holds for the considered model. This result seems to clarify and formalize the general strategy leading to metastability in the pathwise sense followed in all works reviewed in Section \ref{sec:metastability_GLmodel}. It is natural to wonder whether such a result could be used to prove metastability under general assumptions on the different degrees of freedom of the GL model, without having to build an ad hoc proof for each specific instance of the model.

\newpage

\section*{Acknowledgments}
This article is dedicated to the memory of both Christophe Pouzat and Antonio Galves. The latter introduced the two authors to the subject of metastability, and the former was the one who suggested the idea for the present article; he would have participated in its writing had it not been for his unexpected passing.

This work was produced as part of the activities of FAPESP Research, Innovation and Dissemination Center for Neuromathematics (Grant 2013/07699-0) and K. Laxa was supported by a FAPESP fellowship (Grant 2022/07386-0).

\printbibliography

\end{document}